\documentclass[11pt]{article}

\usepackage{amsfonts,amsthm}
\usepackage{mathrsfs}
\usepackage{amsmath,mathtools}   
\usepackage{latexsym}
\usepackage{amssymb}
\usepackage{mathrsfs}
\usepackage{pifont}
\usepackage{tikz,ulem}
\usepackage[all]{xy}
\usepackage{hyperref}
\usepackage{enumitem}

\newtheorem{theo}{Theorem}[section]

\newtheorem*{theo*}{Theorem}
\newtheorem{lemm}[theo]{Lemma}
\newtheorem{prop}[theo]{Proposition}

\newtheorem{conj}[theo]{Conjecture}

\newtheorem{defi}[theo]{Definition}

\newtheorem{remax}[theo]{Remark}
\newenvironment{rema}
  {\pushQED{\qed}\remax}
  {\popQED\endremax}
  
\theoremstyle{definition}

\newtheorem*{term*}{Notation/Terminology}

\newcommand{\bZ}{\mathbb{Z}}

\newcommand{\cD}{\mathcal{D}}

\newcommand{\adj}{\text{adj}}

\newcommand{\un}{1\!\!1}

\usepackage{geometry}
\newcommand{\doi}[1]{\href{https://doi.org/#1}{doi:#1}}
\numberwithin{equation}{section}

\begin{document}

\title{\bf Factorized $A_2$-Leonard pair}

\author{
Nicolas Cramp\'e\textsuperscript{$1$}\footnote{E-mail: crampe1977@gmail.com}~,
Meri Zaimi\textsuperscript{$2$}\footnote{E-mail: meri.zaimi@umontreal.ca}~,
\\[.9em]
\textsuperscript{$1$}
\small Institut Denis-Poisson CNRS/UMR 7013 - Universit\'e de Tours - Universit\'e
d'Orl\'eans,\\
\small~Parc de Grandmont, 37200 Tours, France.\\[.9em]
\textsuperscript{$1$}
\small Laboratoire d'Annecy-le-Vieux de Physique Th\'eorique LAPTh,\\
\small~Universit\'e Savoie Mont Blanc, CNRS, F-74000 Annecy,
 France.\\[.9em]
\textsuperscript{$2$}
\small Centre de Recherches Math\'ematiques, Universit\'e de Montr\'eal, P.O. Box 6128, \\
\small Centre-ville Station, Montr\'eal (Qu\'ebec), H3C 3J7, Canada.\\[.9em]
}
\date{}
\maketitle

\bigskip\bigskip 

\begin{center}
\begin{minipage}{14cm}
\begin{center}
{\bf Abstract}\\
\end{center}
    The notion of factorized $A_2$-Leonard pair is introduced. It is defined as a rank~2 Leonard pair, with actions in certain bases corresponding to the root system of the Weyl group $A_2$, and with some additional properties. The functions arising as entries of transition matrices are bivariate orthogonal polynomials (of Tratnik type) with bispectral properties. Examples of factorized $A_2$-Leonard pairs are constructed using classical Leonard pairs associated to families of orthogonal polynomials of the \mbox{($q$-)Askey} scheme. The most general examples are associated to an intricate product of univariate ($q$-)Hahn and dual ($q$-)Hahn polynomials.
\end{minipage}
\end{center}

\medskip

\begin{center}
\begin{minipage}{13cm}

\textbf{MSC2020 database:} 33C80; 33C45; 16G60
\end{minipage}
\end{center}

\clearpage
\newpage

\section{Introduction}

The purpose of this paper is to study a generalized notion of Leonard pair which involves bivariate orthogonal polynomials.

A Leonard pair is formed by two finite-dimensional linear transformations such that there exists a basis where one is diagonal and the other one is irreducible tridiagonal, and vice versa. This notion has been introduced in \cite{Ter}, and it is closely related to the families of orthogonal polynomials in the terminating branch of the ($q$-)Askey scheme \cite{Koek}. Indeed, the matrix entries when performing a change of basis for a Leonard pair can be expressed in terms of polynomials satisfying both a three-term recurrence relation and a three-term ($q$-)difference equation. Leonard pairs are classified according to these families of orthogonal polynomials, the most general being the $q$-Racah polynomials, which are the discrete analogs of the Askey--Wilson polynomials. The recurrence and $q$-difference operators associated to both families of orthogonal polynomials satisfy the Askey--Wilson algebra \cite{Zhe91}. This algebra, first introduced in \cite{GrZh} to study the symmetry of the $6j$-symbol for the quantum group $U_q(sl_2)$, plays an important role in various related contexts such as centralizers of $U_q(\mathfrak{sl}_2)$ \cite{CVZ,CPVZ}, the Skein algebra \cite{cooke}, or symmetries of physical models \cite{GWZ13,KKM,Post} (for a review see \cite{avatar}). Leonard pairs are recovered when considering finite-dimensional irreducible representations of the Askey--Wilson algebra. The classification of Leonard pairs in terms of orthogonal polynomials of the ($q$-)Askey scheme is an algebraic version of Leonard's theorem \cite{Leo}. This theorem is particularly useful for classifying association schemes with both $P$- and $Q$-polynomial properties in the context of algebraic combinatorics \cite{BI}. The connection between such association schemes with polynomial structures and Leonard pairs is obtained when studying the representations of the subconstituent (also called Terwilliger) algebra \cite{Ter1,Ter2,Ter3}, which is the algebra formed by the adjacency and dual adjacency matrices. 

Multivariate generalizations of the orthogonal polynomials of the ($q$-)Askey scheme have been considered in past works \cite{Griff,Tra,Rosengren,HR,GI,GI2,IX}. Algebras associated to them have also been investigated, for instance the higher rank Bannai--Ito, Racah and Askey--Wilson algebras \cite{BCV,DGVV1,PW,Dec,cooke,icosi,CFPR}. More recently, multivariate generalizations of the $P$- and $Q$-polynomial properties for association schemes have been introduced and studied \cite{BCPVZ,BKZZ,CVZZ,BCVZZ}. However, an analog of Leonard's theorem for multivariate polynomials is not known at the present time. In order to better understand the structure of higher rank association schemes and of multivariate orthogonal polynomials with bispectral properties, and with the purpose of eventually arriving at some classification of these objects, it is crucial to understand the higher rank generalizations of Leonard pairs. 

A generalized notion of Leonard pair has been proposed in \cite{IT}, based on an algebraic interpretation of the Rahman polynomials. Instead of considering a pair of linear transformations with tridiagonal actions in some eigenbases, one considers a pair of $M$-dimensional subspaces of commuting linear transformations with adjacent actions in some common eigenbases. Such pairs are related to $M$-variate polynomials which satisfy recurrence relations with a number of terms that can be higher than three. However, a systematic and complete study of all pairs satisfying the proposed definition seems difficult. In \cite{BCPVZ,CVZZ}, the bivariate polynomial structures of association schemes are characterized by the types of recurrence relations that the associated polynomials satisfy. Some of the natural examples of association schemes examined in \cite{BCPVZ} have a type such that one of the recurrence relations has at most three terms. This is the case for the non-binary Johnson schemes. It is shown in \cite{CVZZ} that the Terwilliger algebra of these schemes satisfies some interesting properties that can be viewed as a refinement of the definition of higher rank Leonard pairs proposed in \cite{IT}. The goal of this paper is to abstract these properties and provide a framework for studying this class of refined higher rank Leonard pairs, that we call factorized Leonard pairs.

The paper will unfold as follows. Section \ref{sec:clasLP} recalls the definition of the classical Leonard pairs introduced in \cite{Ter}, and their description in terms of a parameter array. In Section \ref{sec:intro}, the factorized Leonard pairs are defined as a refinement of the higher rank Leonard pairs of \cite{IT}. We focus in this paper on the case where the adjacency of the actions is of type $A_2$, with reference to the Weyl group root system. Some general properties of these factorized Leonard pairs are also investigated. Section \ref{sec:condtruction} explains how to construct factorized Leonard pairs using classical Leonard pairs. The bivariate polynomials associated to a factorized Leonard pair are discussed in Section \ref{sec:pol}. Examples of factorized $A_2$-Leonard pairs associated to families of orthogonal polynomials of the Askey scheme (type II) and its $q$-analog (type I) are provided respectively in Sections \ref{sec:listA2II} and \ref{sec:listA2I}. One example (associated to the dual ($q$-)Hahn and  ($q$-)Hahn polynomials) in each section is described in detail, the others can be obtained through specializations and limits. In Section \ref{sec:introM}, a higher rank generalized notion of factorized $A_M$-Leonard pair is briefly introduced. Section \ref{sec:conclu} discusses some outlooks. Finally, in Appendix \ref{app:pa}, a list of parameter arrays for the classical Leonard pairs is recalled, and useful contiguity relations for the associated polynomials are provided.

\section{Leonard pairs\label{sec:clasLP}}

We recall in this section the definition of Leonard pairs and their classification \cite{Ter} (see also \cite{TerNote}).
\begin{defi}
Let $V$ denote a vector space over a field $K$ with finite positive dimension. By a Leonard pair on $V$ we mean an ordered pair $(A, A^\star) $, where $A:V\rightarrow V$ and $A^\star:V \rightarrow V$ are linear transformations which satisfy the following two properties:
\begin{itemize}
    \item[(i)] There exists a basis for $V$ with respect to which the matrix representing $A$ is irreducible tridiagonal and the matrix representing $A^\star$ is diagonal;
     \item[(ii)] there exists a basis for $V$ with respect to which the matrix representing $A^\star$ is irreducible tridiagonal and the matrix representing $A$ is diagonal.
\end{itemize}
\end{defi}
In this paper, we restrict ourselves to the case where the field $K$ is of characteristic zero.

The classification of the Leonard pairs is well-known and is  closely related to 
the orthogonal polynomials of the ($q$-)Askey scheme \cite{Ter,TerClas} and to the Askey--Wilson algebra \cite{TV03}. 
We follow here the notations in \cite{TerClas}.

To each Leonard pair $(A,A^\star)$ with $\dim(V)=d+1$ ($d$ is called the diameter of the Leonard pair), a
set of parameters is associated  
\begin{align}\label{eq:Phi}
\Phi=(\theta_i,\theta^\star_i,i=0,\dots,d;\varphi_j, \phi_j,j=1,\dots,d)\,,
\end{align}
called parameter array. To be precise, there is a one-to-one correspondence between the parameter arrays and the Leonard systems, which are Leonard pairs with orderings of the primitive idempotents. Each Leonard pair is associated to four Leonard systems. We shall always assume that an ordering of the primitive idempotents is fixed for a given Leonard pair such that it forms a Leonard system, but for simplicity we will not mention explicitly the Leonard system.

From the parameter array \eqref{eq:Phi}, one gets:
\begin{itemize}
\item $\{ \theta_i\ |\ i=0,\dots, d \}$ are the eigenvalues of $A$;
\item $\{ \theta^\star_i\ |\ i=0,\dots, d \}$ are the eigenvalues of $A^\star$;
\item in the basis $\{ v_0^\star,\dots, v_d^\star\}$ where $A$ is tridiagonal and $A^\star$ is diagonal, they are written as 
\begin{align}\label{eq:Atri}
&A v_i^\star= b_{i-1} v_{i-1}^\star +a_{i} v_{i}^\star+ c_{i+1} v_{i+1}^\star \, ,\\
&A^\star v_i^\star = \theta_i^\star v_i^\star\,,
\end{align}
where, by convention $b_{-1}=c_{d+1}=0$ and 
\begin{subequations}\label{eq:abc}
\begin{align}
& a_{i}=a_i(\Phi)=\theta_i +\frac{\varphi_i}{\theta_i^\star-\theta_{i-1}^\star}
+\frac{\varphi_{i+1}}{\theta_i^\star-\theta_{i+1}^\star}\,,\\
&b_i=b_i(\Phi)= \varphi_{i+1} \frac{\displaystyle \prod_{h=0}^{i-1}(\theta^\star_{i}-\theta^\star_{h})} {\displaystyle  \prod_{h=0}^{i}(\theta^\star_{i+1}-\theta^\star_{h} ) }\,,
\qquad c_i= c_i(\Phi)= \phi_{i}\frac{\displaystyle \prod_{h=0}^{d-i-1}(\theta^\star_{i}-\theta^\star_{d-h})} {\displaystyle  \prod_{h=0}^{d-i}(\theta^\star_{i-1}-\theta^\star_{d-h} ) }\,;
\end{align}
\end{subequations}
\item in the basis  $\{ v_0,\dots, v_d\}$  where $A$ is diagonal and $A^\star$ is tridiagonal, they are written as 
\begin{align}
&A v_x= \theta_x v_x \, ,\label{eq:A1}\\
&A^\star v_x = b^\star_{x-1} v_{x-1} +a^\star_{x} v_{x}+ c^\star_{x+1} v_{x+1} \,,\label{eq:A2}
\end{align}
where, by convention $b^\star_{-1}=c^\star_{d+1}=0$, 
\begin{align}\label{eq:abc2}
& a^\star_{i}=a_{i}(\Phi^\star)\,,\qquad b^\star_i=b_{i}(\Phi^\star)\,,
\qquad c^\star_i= c_{i}(\Phi^\star)\,,
\end{align}
and the transformed parameter array is
 \begin{align}\label{eq:Phistar}
     \Phi^\star=(\theta^\star_i,\theta_i,i=0,\dots,d;\varphi_j, \phi_{d-j+1},j=1,\dots,d)\,.
 \end{align}
\end{itemize}
The classification of Leonard pairs is obtained by providing the complete list of all the possible parameter arrays. We recall the parameter arrays used in this paper in Appendix \ref{app:pa}.

\begin{rema}\label{rem:norma} Instead of relations \eqref{eq:A1}-\eqref{eq:A2}, one can have 
\begin{align}
&A v_x= \theta_x v_x \, ,\\
&A^\star v_x =\frac{\alpha_{x-1}}{\alpha_x} b^\star_{x-1} v_{x-1} +a^\star_{x} v_{x}+ \frac{\alpha_{x+1}}{\alpha_x} c^\star_{x+1} v_{x+1} ,
\end{align}
for $\alpha_x$ non-vanishing free parameters. They correspond to the normalizations of the vectors $v_x$. Similarly one can modify the normalizations of $v^\star_i$ such that relation \eqref{eq:Atri} is modified. 
\end{rema}

For any Leonard pair, the change of basis between both bases is known explicitly and is given in terms of orthogonal polynomials of the ($q$-)Askey scheme.
Explicitly, to each parameter array $\Phi$, we associate the following polynomial of degree $i$ w.r.t. $\lambda$ \cite{TerClas}:
\begin{align}
    u_i(\lambda)=\sum_{k=0}^i \prod_{j=0}^{k-1}
    \frac{(\lambda-\theta_j) (\theta^\star_i-\theta^\star_j)}
    {\varphi_{j+1}}\,.
\end{align}
Let us define
\begin{align}\label{eq:pxi}
  p_{xi}(\Phi)= k_i(\Phi) u_i(\theta_x)\,,
\end{align}
where
\begin{align} \label{eq:kj}
    k_i(\Phi)=\prod_{h=1}^i\frac{\varphi_{h}}{\phi_{h}}\ \frac{\displaystyle  \prod_{h=1}^{d}(\theta_0^\star-\theta_{h}^\star) }
    {\displaystyle  \prod_{\genfrac{}{}{0pt}{}{h=0}{h\neq i}}^{d}(\theta_i^\star-\theta_{h}^\star)}\, .
\end{align}
Let us remark that $p_{xi}(\Phi)$ is denoted by $v_i(\theta_x)$ in \cite{TerClas}.
These coefficients $p_{xi}(\Phi)$ allow us to connect the two bases:
\begin{equation}
    v_i^\star= \sum_{x=0}^d p_{xi}(\Phi) v_x\,.
 \end{equation}
These coefficients are given explicitly in Appendix \ref{app:pa} for different parameter arrays, and satisfy the recurrence and the difference relations:
\begin{subequations}
\begin{align}
\theta_x\ p_{xi}(\Phi)&= b_{i-1}p_{x,i-1}(\Phi)
+a_{i}p_{xi}(\Phi)+c_{i+1}p_{x,i+1}(\Phi)\,, \label{eq:recu}\\
\theta_i^\star\ p_{xi}(\Phi)&=b_x^\star  p_{x+1,i}(\Phi)+a_x^\star  p_{xi}(\Phi)+c_x^\star  p_{x-1,i}(\Phi)\,.\label{eq:diff}
\end{align}
\end{subequations}
They also satisfy the following orthogonality relation
\begin{align}
\sum_{i=0}^d p_{xi}(\Phi) p_{iy}(\Phi^\star) =\delta_{xy}\, \nu(\Phi)\,, \label{eq:Pinv}
\end{align}
where $\Phi^\star$ is given by \eqref{eq:Phistar}  and 
\begin{align}
&  \nu(\Phi)=\nu(\Phi^\star)=\prod_{h=1}^{d}\frac{(\theta_0-\theta_{h})(\theta^\star_0-\theta^\star_{h})}{\phi_{h}}\,.\label{eq:nu}
\end{align}

\section{Definitions and properties}\label{sec:intro}

In this section, we provide the definition of a factorized Leonard pair, generalizing the notion of Leonard pairs as well as some of its properties.

\subsection{Definition of a factorized $A_2$-Leonard pair }\label{sec:intro1}

Let $N$ be a positive integer and $\cD$ be the set constituted of couples of non-negative integers whose sum is at most $N$:
\begin{equation}\label{eq:dom}
 \cD=\{ (x,y)\ | \ x,y\in \bZ_{\geq 0}, \ x+y\leq N \}\,.
\end{equation}
Two couples $(x,y),(x',y')\in \cD$ are called $A_2$-adjacent, $(x,y)\ \adj\ (x',y')$, if
\begin{equation}
  (x-x',y-y') \in A_2,
\end{equation}
where $A_2$ is the following subset of $\bZ^2$:
\begin{align} \label{eq:AB2}
 &A_2=\{(0,0),(1,-1), (1,0), (0,1),  (-1,1), (-1,0), (0,-1)\}\,.
 \end{align}
The notation $A_2$ for the previous set comes from the fact that it corresponds to the root system of the Weyl group $A_2$.

\begin{defi}\label{def:IT}\cite{IT}
Let $V$ be a $\mathbb{C}$-vector space.
The pair $(H,H^\star)$ is a $A_2$-Leonard pair on the domain $\cD$ if the following statements are satisfied:
\begin{itemize}
 \item[(i)] $H$ is a two-dimensional  subspace of $\text{End}(V)$ whose elements are diagonalizable and mutually commute;
 \item[(ii)] $H^\star$ is a two-dimensional  subspace of $\text{End}(V)$ whose elements are diagonalizable and mutually commute;
 \item[(iii)] there exists a bijection $\lambda\mapsto V_\lambda$  from $\cD$ to the common eigenspaces of $H$ such that, for all $\lambda\in \cD$,
 \begin{equation}\label{eq:HsV}
  H^\star V_\lambda \subseteq \sum_{ \genfrac{}{}{0pt}{}{\chi \in \cD}{\chi\ \adj \  \lambda}} V_\chi\,;
 \end{equation}
 \item[(iv)]  there exists a bijection $\lambda\mapsto V^\star_\lambda$  from $\cD$ to the common eigenspaces $V^\star_\lambda$ of $H^\star$ such that, for all $\lambda\in \cD$,
 \begin{equation}\label{eq:HVs}
  H V^\star_\lambda \subseteq \sum_{ \genfrac{}{}{0pt}{}{\chi \in \cD}{\chi\ \adj \  \lambda}} V^\star_\chi\,;
 \end{equation}
 \item[(v)] there does not exist a subspace $W$ of $V$ such that $HW\subseteq W$, $H^\star W\subseteq W$, $W\neq 0$, $W\neq V$;
 \item[(vi)] $\text{dim}(V_\lambda)=\text{dim}(V^\star_\lambda)=1$, for all  $\lambda\in \cD$;
 \item[(vii)] there exists a non-degenerate symmetric bilinear form  $\langle,\rangle$ on $V$ such that, for $\lambda,\chi\in \cD$ and $\lambda\neq \chi$,
 \begin{equation}
\langle V_\lambda, V_\chi\rangle=0\ ,\qquad  \langle V^\star_\lambda, V^\star_\chi\rangle=0\,.
 \end{equation}
\end{itemize}
\end{defi}
Let us notice that the dimension of the vector space $V$ is finite and is given more precisely by $\dim(V) = |\cD|=\binom{N+2}{2}$. 
\begin{rema} The previous definition has been given in the conclusion of \cite{IT}. The definition of adjacency used in this paper is slightly modified in comparison to \cite{IT}, which explains the difference in relations \eqref{eq:HsV} and \eqref{eq:HVs}.
\end{rema}

The classification of the usual Leonard pairs has been recalled in Section \ref{sec:clasLP} and is closely related to 
the orthogonal polynomials of the ($q$-)Askey scheme.
The classification of the $A_2$-Leonard pairs would provide similar results but for bivariate polynomials. For a first step in this direction, we propose to add further conditions to the previous definition such that a classification seems possible.
We arrive at the following definition of factorized Leonard pairs.

\begin{defi} \label{def:flp}
Let $(H,H^\star)$ be a $A_2$-Leonard pair on the domain $\cD$ and $X,Y\in \mathrm{End}(V)$ (resp. $X^\star,Y^\star$) be a basis of $H$ (resp. $H^\star$).
The pair $(H,H^\star)$ is a \textit{factorized $A_2$-Leonard pair} if the following statements are satisfied:
\begin{itemize}
 \item[(viii)] the eigenvalues of $X$ and $Y^\star$ satisfy, for any $(x,y),(i,j)\in \cD$,
 \begin{align}
 & (X-t_x \un) V_{(x,y)} =0 \,,  \qquad 
 (Y^\star-\theta^\star_j \un) V^\star_{(i,j)} =0 \,,
  \end{align}
  where $\un$ is the identity in $\mathrm{End}(V)$ and $t_x\neq t_{x'} \Leftrightarrow x\neq x'$ 
  and $\theta^\star_j\neq \theta^\star_{j'} \Leftrightarrow j\neq j'$;
 \item[(ix)] in addition to the commutation relations $[X,Y]=0$ and $[X^\star,Y^\star]=0$, one has
 \begin{align}\label{eq:com1}
  [X,Y^\star]=0\,;
 \end{align}
 \item[(x)] each common eigenspace of $X$ and $Y^\star$ has dimension one;
   \item[(xi)]   
  the following relations hold, for $(i,j),(i+1,j)\in \cD$ and $(x,y),(x,y+1)\in \cD$,
 \begin{subequations}
      \begin{align}
 & \langle\, V^\star_{(i+1,j)}\, ,\, X V^\star_{(i,j)}\,\rangle\neq 0\,,&& \langle\, V^\star_{(i,j)}\, ,\, X V^\star_{(i+1,j)}\,\rangle\neq 0\,, \label{eq:irre1}\\
      &\langle V_{(x,y+1)}\,,\, Y^\star V_{(x,y)}\rangle \neq 0\,,&& \langle V_{(x,y)}\,,\, Y^\star V_{(x,y+1)}\rangle \neq 0\,.\label{eq:irre2}
  \end{align} 
  \end{subequations}
 \end{itemize}
 The ordered 4-tuple $(X,Y;X^\star,Y^\star)$ is called a basis of the factorized $A_2$-Leonard pair.
 \end{defi}

\begin{rema}
In the previous definition, one could replace the pair $(X,Y^\star)$ by $(X^\star,Y)$. 
\end{rema}

In the following, we shall use the shortened name factorized Leonard pair instead of factorized $A_2$-Leonard pair.

\subsection{Properties of a factorized Leonard pair }\label{sec:intro2}

In this section, let $(H,H^\star)$ be a factorized Leonard pair with basis $(X,Y;X^\star,Y^\star)$. We provide different properties for this pair, allowing us to study it in detail.

The first property concerns a symmetry of its basis.
\begin{prop}\label{prop:cb}
 For $a_1,a_2,a^\star_1,a^\star_2$ invertible elements of the field $K$ and $b_1,b_2,b_1^\star,b_2^\star, c,c^\star\in K$, the 4-tuple
    \begin{equation}
    (a_1 X+b_1 \un ,\ a_2 Y+c X+b_2\un\ ;\ a^\star_1 X^\star +c^\star Y^\star +b^\star_1\un ,\ a_2^\star Y^\star+b_2^\star \un)\,,\label{eq:nb1}
    \end{equation}
    where $\un$ is the identity in $End(V)$, 
    is another basis of the factorized Leonard pair $(H,H^\star)$.
\end{prop}
\proof It is straightforward to check that all the properties of the definition of a factorized Leonard pair are satisfied by \eqref{eq:nb1}.   \endproof

To study further the factorized Leonard pairs, it is useful to introduce some projectors. 
Indeed, let us associate to the factorized Leonard pair the following projectors:
\begin{align}\label{eq:projectors}
 E_x= \prod_{\genfrac{}{}{0pt}{}{x'=0}{x'\neq x}  }^N  \frac{X-t_{x'}}{t_x-t_{x'}}\,, \qquad E^\star_j= \prod_{\genfrac{}{}{0pt}{}{j'=0}{j'\neq j} }^N  \frac{Y^\star-\theta^\star_{j'}}{\theta^\star_j-\theta^\star_{j'}}\,,
\end{align}
which are such that $XE_x=t_xE_x$ and $Y^\star E^\star_j=\theta^\star_j E_j^\star$.
It follows from \eqref{eq:com1} and \eqref{eq:projectors} that
\begin{equation}
E_xE_{x'}=\delta_{xx'}E_x\,,\qquad E_j^\star E_{j'}^\star=\delta_{jj'}E_j\,,\qquad  \sum_{x=0}^NE_x=\un \,,\qquad  \sum_{j=0}^N E^\star_j=\un \,,\qquad [E_x , E^\star_j]=0\,.
\end{equation}

\begin{lemm}\label{lem:dimtildeV}
For $0\leq x,j\leq N$, the dimension of the vector space $E_x E^\star_jV$ is given by
\begin{align}
\dim(E_x E^\star_jV) =
\begin{cases}
1\,, & \text{if } x+j\leq N,\\ 0\,, & \text{otherwise}.   
\end{cases}
\end{align}
\end{lemm}
\proof $E_x E^\star_jV$ is a common vector space of $X$ and $Y^\star$ or is zero. Then due to the point (x) of Definition \ref{def:flp}, one deduces that, for $0\leq x,j\leq N$,
\begin{align}\label{eq:dimExEsj}
&\dim( E_x E^\star_jV )\leq 1\,.
\end{align}
Let us remark that $\dim(E_0V)=N+1$ (due to the point (iii) of Definition \ref{def:IT}).
It follows from this remark and equation \eqref{eq:dimExEsj} that for any $j=0,1,\dots, N$,
\begin{align}\label{eq:dimE0Esj}
\dim(E_0 E^\star_jV)=1\,.
\end{align}
Since $\dim(E_1V)=N$ and $\dim(E^\star_N V)=1$ (due to the points (iii), (iv) of Definition \ref{def:IT}), one deduces using equations \eqref{eq:dimExEsj} and \eqref{eq:dimE0Esj} that for any $j=0,1,\dots, N-1$,
\begin{align}
\dim( E_1 E^\star_jV)=1\,, \quad \text{and} \qquad \dim( E_1 E^\star_N V)=0\,.
\end{align}
Recursively, one demonstrates the result of the lemma.
\endproof

From the previous lemma, one deduces that each common eigenspace $\widetilde V_{(x,j)}$ ($0\leq x,j\leq N$ and $x+j\leq N$) of $X$ and $Y^\star$ can be written as follows:
\begin{align}
    \widetilde V_{(x,j)} = E_x E^\star_jV = E^\star_jE_x V\,,
\end{align}
such that
\begin{align}
 &\big( X-t_x \un \big) \widetilde V_{(x,j)}= 0 \,,\qquad
 \big( Y^\star-\theta^\star_j \un \big)  \widetilde V_{(x,j)}= 0\,.
\end{align}

\begin{lemm} \label{lem:LP}
For any $j=0,1,\dots,N$, $(X,X^\star)$ restricted on $E_j^\star V$ is a Leonard pair with diameter $N-j$ and, for any $x=0,1,\dots,N$, $(Y,Y^\star)$ restricted on $E_x V$ is a Leonard pair with diameter $N-x$.
 \end{lemm}
\proof 
Let us fix $j \in \{0,1,\dots,N\}$. On one hand, the subspace $E_j^\star V$ can be written as
\begin{equation}
E_j^\star V = \sum_{\lambda \in \cD} E_j^\star V_\lambda^\star = \sum_{i=0}^{N-j}  V_{(i,j)}^\star \, .
\end{equation}
The subspace $V_{(i,j)}^\star$ is invariant under $X^\star$ for $(i,j)\in\cD$.
Using condition (iv) of Definition \ref{def:IT} and the commutation relation $[X,Y^\star]=0$, one finds for any $i=0,1,\dots,N-j$
\begin{equation}
X V_{(i,j)}^\star = E_j^\star X V_{(i,j)}^\star \subseteq  E_j^\star \sum_{\lambda \  \adj \ (i,j)} V_\lambda^\star = V_{(i,j)}^\star +  V_{(i+1,j)}^\star +  V_{(i-1,j)}^\star \,. 
\end{equation}
Recall that the spaces $V_{(i,j)}^\star$ have dimension $1$ from condition (vi) of Definition \ref{def:IT} and denote $V_{i,j}^\star$ one non-vanishing vector of this vector space.
The set $\{V_{i,j}^\star\}_{i=0}^{N-j}$ is a basis of $E_j^\star V$. In this basis, $X$ is represented by a tridiagonal matrix and $X^\star$ by a diagonal matrix. 
In addition, the matrix representing $X$ is an irreducible tridiagonal matrix due to relations \eqref{eq:irre1} \textit{i.e.} there are no vanishing entries in the lower and upper diagonals. One deduces that there does not exist a subspace $W_j$ of $E_j^\star V$ such that $XW_j\subseteq W_j$, $X^\star W_j\subseteq W_j$, $W_j\neq 0$, $W_j\neq E_j^\star V$ (the fact that all the eigenvalues of $X^\star$ restricted to $E_j^\star V$ are different two by two has also been used).

On the other hand, one can write 
\begin{equation}
E_j^\star V = \sum_{\lambda \in \cD} E_j^\star \widetilde V_\lambda = \sum_{x=0}^{N-j}  \widetilde V_{(x,j)} \, . 
\end{equation}
The subspace $\widetilde V_{(x,j)}$ is invariant under $X$ and one can use condition (iii) of Definition \ref{def:IT} together with $[X^\star,Y^\star]=0$ to get for any $x=0,1,\dots,N-j$
\begin{equation}
X^\star \widetilde V_{(x,j)} = E_j^\star X^\star E_x V \subseteq E_j^\star (E_{x} V + E_{x+1} V + E_{x-1} V)   = \widetilde V_{(x,j)}+\widetilde V_{(x+1,j)}+\widetilde V_{(x-1,j)} \, .
\end{equation}
Let us denote $\widetilde V_{x,j}$ a basis vector of the one-dimensional vector space $\widetilde V_{(x,j)}$. The set $\{\widetilde V_{x,j}\}_{x=0}^{N-j}$ is a basis of $E_j^\star V$.
In this basis,  $X$ is represented by a diagonal matrix and $X^\star$ by an irreducible tridiagonal matrix. It follows that $X$ and $X^\star$ restricted on $E_j^\star V$ form a Leonard pair. 

The proof that $(Y,Y^\star)$ restricted on $E_x V$ is a Leonard pair is similar. In this case, $Y$ is diagonal and $Y^\star$ is irreducible tridiagonal on $\{ V_{(x,y)}\}_{y=0}^{N-x}$, and $Y$ is irreducible tridiagonal and $Y^\star$ is diagonal on $\{\widetilde V_{(x,j)}\}_{j=0}^{N-x}$.
\endproof

\begin{lemm}
   For a factorized Leonard pair, statement $(v)$ of Definition \ref{def:IT} concerning the irreducibility is redundant.
\end{lemm}
\proof The subspace $\widetilde V_{(x,j)}$ is invariant under $X$ and $Y^\star$, and the couples of eigenvalues $(t_x,\theta_j^\star)$ are two by two different. Moreover, from the proof of the previous lemma, one gets, for $(x,j),(x+1,j)\in \cD$,
\begin{subequations}
     \begin{align}
 & \langle \widetilde V_{(x+1,j)} , X^\star \widetilde V_{(x,j)}\rangle \neq 0\,,&& \langle \widetilde V_{(x,j)}, X^\star \widetilde V_{(x+1,j)}\rangle\neq 0\,, \label{eq:irre7}
 \end{align}
 and, for $(x,j),(x,j+1)\in \cD$,
 \begin{align}
      &\langle \widetilde V_{(x,j+1)} ,Y \widetilde V_{(x,j)}\rangle\neq 0\,,&&\langle \widetilde V_{(x,j)}, Y \widetilde V_{(x,j+1)}\rangle\neq 0\,.\label{eq:irre8}
  \end{align}
\end{subequations}
Therefore one deduces that $X$, $Y$, $X^\star$ and $Y^\star$ act irreducibly on $V$.
\endproof 

\section{Construction of factorized Leonard pairs \label{sec:condtruction}}

In this section, we describe a general method for constructing factorized Leonard pairs.

\subsection{Parameter arrays associated to a factorized Leonard pair}

\begin{figure}[htbp]
\begin{center}
\begin{minipage}[t]{.8\linewidth}
\begin{center}
\begin{tikzpicture}[scale=1]
\draw[-] (-5,0)--(0,0)--(5,0);
\draw [fill] (-5,0) circle (0.1);
\draw (-5,0) node[above] {$X,Y$};
\draw (-5,0) node[below] {$V_{x,y}$};

\draw[-] (-2.6,0.1)--(-2.5,0)--(-2.6,-0.1);
\draw (-2.5,0) node[above] {$\Psi_{x}$};
\draw (-2.5,0) node[below] {$p_{yj}(\Psi_x)$};

\draw [fill] (0,0) circle (0.1); 
\draw (0,0) node[above] {$Y^\star,X$};
\draw (0,0) node[below] {$\widetilde V_{x,j}$};

\draw[-] (2.4,0.1)--(2.5,0)--(2.4,-0.1);
\draw (2.5,0) node[above] {$\Phi_{j}$};
\draw (2.5,0) node[below] {$p_{xi}(\Phi_j)$};

\draw [fill] (5,0) circle (0.1);  
\draw (5,0) node[above] {$X^\star,Y^\star$};
\draw (5,0) node[below] {$V^\star_{i,j}$};

\end{tikzpicture}
\caption{Above each vertex, a couple of commuting operators is written and, below, the vectors of the basis where these operators are diagonal. At each oriented edge is associated the parameter array of the Leonard pair given by the operators above this edge.  \label{fig:tr}}
\end{center}
\end{minipage}
\end{center}
\end{figure}
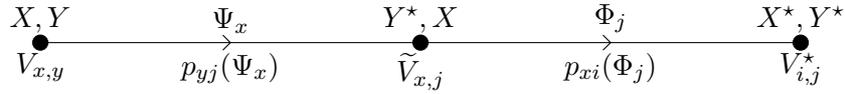

Let $(X,Y;X^\star,Y^\star)$ be a basis of a factorized Leonard pair. 
From Lemma \ref{lem:LP}, $(X,X^\star)$ restricted on $E^\star_j V$ is a Leonard pair of diameter $N-j$, for any $j=0,\dots,N$.
From the result recalled in Section \ref{sec:clasLP}, we can associate to each of these Leonard pairs a parameter array given explicitly by
\begin{subequations}\label{eq:paflp}
\begin{align} \label{eq:pa1}
&\Phi_j=(t_x ,t^\star_x(j),x=0,\dots,N-j;\ f_x(j), g_x(j),x=1,\dots,N-j)\,.
\end{align}
Similarly $(Y,Y^\star)$ restricted on $E_x V$ is a Leonard pair of diameter $N-x$, for any $x=0,\dots,N$ and we can associate to each of these pairs the following parameter array
\begin{align} 
\label{eq:pa2}
&\Psi_x=(\theta_j(x) ,\theta^\star_j,j=0,\dots,N-x;\ \varphi_j(x), \phi_j(x),j=1,\dots,N-x)\,.
\end{align}
\end{subequations}
Let us emphasize that, for a factorized Leonard pair, the parameters in $\Phi_j$ can depend explicitly on $j$ except the eigenvalues $t_x$ because of statement (viii) of Definition \ref{def:flp}. Similarly, the dual eigenvalues in $\Psi_x$ do not depend on $x$. 
We can visualize the different parameter arrays associated to the different Leonard pairs as in Figure \ref{fig:tr}.

Due to Proposition \ref{prop:cb}, the following transformations on the parameter arrays \eqref{eq:paflp} associated to a factorized Leonard pair, 
\begin{subequations}\label{eq:inv}
\begin{align}
&\Phi_j=(a_1 t_x+b_1 ,a_2 t^\star_x(j)+c t_x+b_2,x=0,\dots,N-j;\ a_1 a_2 f_x(j), a_1 a_2 g_x(j),x=1\dots N-j)\,, \\
&\Psi_x=(a_1^\star \theta_j(x)+ c^\star \theta^\star_j+b_1^\star,a_2^\star\theta^\star_j+b_2^\star,j=0\dots N-x;a_1^\star a_2^\star\varphi_j(x), a_1^\star a_2^\star\phi_j(x),j=1\dots N-x)\,,
\end{align}
\end{subequations}
 with $a_1,a_2,a_1^\star,a_2^\star,b_1,b_2,b_1^\star,b_2^\star,c,c^\star \in K$ and $a_1,a_2,a_1^\star,a_2^\star$ invertible, are also the associated parameter arrays of a factorized Leonard pair. \\

From the knowledge of the parameter arrays, one can obtain the actions of the elements of the factorized Leonard pair basis $(X,Y;X^\star,Y^\star)$ on the different bases of the vector space $V$. As before, we denote by $V_{x,y}$ (resp. $\widetilde V_{x,j}$, $V^\star_{i,j}$) a basis vector of the one-dimensional space  $V_{(x,y)}$ (resp. $\widetilde V_{(x,j)}$, $V^\star_{(i,j)}$).
These different vectors are also displayed in Figure \ref{fig:tr}. In the eigenspace $E^\star_jV$, one gets $(Y^\star-\theta^\star_j\un)E^\star_jV=0$ and thus 
\begin{equation}
    Y^\star \widetilde V_{x,j} = \theta^\star_j \widetilde V_{x,j}\,, \qquad 
   Y^\star V^\star_{i,j} = \theta^\star_j V^\star_{i,j}\,.
\end{equation}
By definition of the parameter array $\Phi_j$ one gets:
\begin{align}
    &X \widetilde V_{x,j} = t_x \widetilde V_{x,j}\,,\qquad  X V^\star_{i,j}=b_{i-1}(\Phi_j) V^\star_{i-1,j}+
 a_{i}(\Phi_j) V^\star_{i,j}+c_{i+1}(\Phi_j) V^\star_{i+1,j}\,,\\
& X^\star \widetilde V_{x,j}= 
c_{x+1}(\Phi^{\star}_j)\ \widetilde V_{x+1,j}+a_x(\Phi^{\star}_j)\ \widetilde V_{x,j}+
 b_{x-1}(\Phi^{\star}_j)\ \widetilde V_{x-1,j}\,, 
\qquad X^\star V^\star_{i,j} = t^\star_i(j) V^\star_{i,j}\,.
\end{align}
In the eigenspace $E_x V$, one gets $X E_x V =t_x E_x V$ and thus 
\begin{equation}
    X V_{x,y} = t_x V_{x,y}\,,\qquad X \widetilde V_{x,j} = t_x \widetilde V_{x,j}\,,
\end{equation}
and by definition of the parameter array $\Psi_x$, one gets:
\begin{align}
& Y V_{x,y}= \theta_y(x) V_{x,y}\,,\qquad
Y \widetilde V_{x,j}=c_{j+1}(\Psi_x)\widetilde V_{x,j+1}
+ a_{j}(\Psi_x)\widetilde V_{x,j}
+ b_{j-1}(\Psi_x)\widetilde V_{x,j-1}\,, \label{YVt}\\
   & Y^\star V_{x,y}=b_{y-1}(\Psi_x^\star) V_{x,y-1}+
a_{y}(\Psi_x^\star) V_{x,y}+ c_{y+1}(\Psi_x^\star) V_{x,y+1}\,,\qquad Y^\star \widetilde V_{x,j}=\theta^\star_j\widetilde V_{x,j}\,.
\end{align}

\begin{rema} As explained in Remark \ref{rem:norma}, an element of a Leonard pair can always be conjugated by a diagonal matrix. Therefore the actions of $X^\star$ and $Y$ on $\widetilde V_{x,j}$ can be transformed as follows:
    \begin{align}
& X^\star \widetilde V_{x,j}= \frac{\alpha^\star_{x+1}(j)}{\alpha^\star_{x}(j)}\  c_{x+1}(\Phi^{\star}_j)\ \widetilde V_{x+1,j}+a_x(\Phi^{\star}_j)\ \widetilde V_{x,j}+
\frac{\alpha^\star_{x-1}(j)}{\alpha^\star_{x}(j)}\ b_{x-1}(\Phi^{\star}_j)\ \widetilde V_{x-1,j}\,,\\
&Y \widetilde V_{x,j}=\frac{\alpha_{j+1}(x)}{\alpha_{j}(x)}\ c_{j+1}(\Psi_x)\ \widetilde V_{x,j+1}
+ a_{j}(\Psi_x)\ \widetilde V_{x,j}
+\frac{\alpha_{j-1}(x)}{\alpha_{j}(x)}\ b_{j-1}(\Psi_x)\ \widetilde V_{x,j-1}\,.
\end{align}
These new coefficients could help us to construct a factorized Leonard pair. However, in all the examples found these coefficients are not useful. To simplify the presentation of this article, we choose to present the results without these coefficients.
\end{rema}

The entries of the transition matrix associated to the Leonard pair $(X,X^\star)$ restricted on $E^\star_j V$ are given by $p_{xi}(\Phi_j)$, defined in \eqref{eq:pxi}. Explicitly, one has 
\begin{equation}\label{eq:Vs2}
V^\star_{i,j}= \sum_{x=0}^{N-j}  p_{xi}(\Phi_j) \widetilde V_{x,j}\,,\qquad \widetilde V_{x,j}=\frac{1}{\nu(\Phi_j)}\sum_{i=0}^{N-j}  p_{ix}(\Phi_j^\star)V^\star_{i,j}\,.
\end{equation}
In the second relation above, we have used the inverse of the transition matrix and relation \eqref{eq:Pinv}. 
The entries $p_{xi}(\Phi_j)$ are indicated below the edge associated to the pair $(X,X^\star)$ in Figure \ref{fig:tr}.
Similarly, $p_{yj}(\Psi_x)$ is the polynomial associated to the parameter array \ref{eq:pa2} such that
\begin{equation}\label{eq:VVt}
\widetilde V_{x,j} =\sum_{y=0}^{N-x} p_{yj}(\Psi_x) V_{x,y}\,,\qquad  V_{x,y}= \frac{1}{\nu(\Psi_x)}\sum_{j=0}^{N-x} p_{jy}(\Psi_x^\star) \widetilde V_{x,j}\,.
\end{equation}
The previous relations allow us to prove the following proposition.
\begin{prop}
The change of basis from $V_{x,y}$ to $V^\star_{i,j}$ is given by
\begin{align}\label{eq:VsV}
V^\star_{i,j}=\sum_{x=0}^{N-j}\sum_{y=0}^{N-x} p_{xi}(\Phi_j)p_{yj}(\Psi_x) V_{x,y}\,.
\end{align}
\end{prop}
The fact that this change of basis associated to a factorized Leonard pair is written as a product of two polynomials justifies the name \textit{factorized}.

It remains to compute the actions of $Y$ on $V^\star_{i,j}$, and of $X^\star$ on $V_{x,y}$.

\subsection{Action of Y on $ V^\star_{i,j}$}

Using the change of basis \eqref{eq:Vs2} and the action \eqref{YVt}, one gets a first expression for the action of Y on $ V^\star_{i,j}$
\begin{align}
Y V^\star_{i,j}&=
\sum_{x=0}^{N-j}p_{xi}(\Phi_j) Y \widetilde V_{x,j}\\
&=\sum_{x=0}^{N-j} p_{xi}(\Phi_j)
\left( c_{j+1}(\Psi_x)\widetilde V_{x,j+1}
+ a_{j}(\Psi_x)\widetilde V_{x,j}
+b_{j-1}(\Psi_x)\widetilde V_{x,j-1}\right)\,.\label{eq:Vt1}
\end{align}
In order to obtain a factorized Leonard pair, this action has to be written as (see statement (iv) of the definition):
\begin{align}
Y V^\star_{i,j}&=\sum_{(\epsilon,\epsilon')\in A_2} C^{\epsilon,\epsilon'}_{i,j}\ V^\star_{i+\epsilon,j+\epsilon'}=\sum_{(\epsilon,\epsilon')\in A_2} C^{\epsilon,\epsilon'}_{i,j}
\sum_{x=0}^{N-j-\epsilon'} p_{x,i+\epsilon}(\Phi_{j+\epsilon'})  \widetilde V_{x,j+\epsilon'} \, ,\label{eq:Vt2}
\end{align}
where all the coefficients $C_{i,j}$ must be independent of $x$ and $y$.
Comparing \eqref{eq:Vt1} and \eqref{eq:Vt2}, one finds
\begin{subequations}\label{eq:cons1}
\begin{align}
   &c_{j+1}(\Psi_x)
   p_{xi}(\Phi_j)
 =\sum_{\epsilon=-1,0} C^{\epsilon,+}_{i,j} p_{x,i+\epsilon}(\Phi_{j+1})\,,\label{eq:cons11}\\
&  a_{j}(\Psi_x) p_{xi}(\Phi_j)
 =\sum_{\epsilon=0,\pm 1} C^{\epsilon,0}_{i,j} p_{x,i+\epsilon}(\Phi_{j})\,,\label{eq:cons12}\\
   &b_{j-1}(\Psi_x)
   p_{xi}(\Phi_j)
 =\sum_{\epsilon=0,+1} C^{\epsilon,-}_{i,j} p_{x,i+\epsilon}(\Phi_{j-1}) \,.\label{eq:cons13}
\end{align}
\end{subequations}
To get a factorized Leonard pair, the previous constraints must be solved. These provide strong restrictions on the parameter arrays which are admissible.

\subsection{Action of $X^\star$ on $ V_{x,y}$}

Similarly as in the previous subsection, the action of $X^\star$ on $ V_{x,y}$ provides constraints between the parameter arrays. Here one must use the change of basis \eqref{eq:VVt}.

Demanding that the action of $X^\star$ on $V_{x,y}$ takes the following form (see statement (iii) of the definition)
\begin{align}
    X^\star V_{x,y} =\sum_{(\epsilon,\epsilon')\in A_2} D^{\epsilon,\epsilon'}_{x,y}V_{x+\epsilon,y+\epsilon'}\,,
\end{align}
where all the coefficients $D_{x,y}$ must be independent of $i$ and $j$, leads to the following constraints
\begin{subequations}\label{eq:cons2}
\begin{align}
   &c_{x+1}(\Phi_j^\star)  p_{jy}(\Psi^\star_x)
= \frac{\nu(\Psi_x)}{\nu(\Psi_{x+1})}\sum_{\epsilon'=-1,0} D^{+,\epsilon'}_{x,y} p_{j,y+\epsilon'}(\Psi^\star_{x+1})\,,\label{eq:cons21}\\
&   a_{x}(\Phi_j^\star) p_{jy}(\Psi^\star_x)
 =\sum_{\epsilon'=0,\pm 1} D^{0,\epsilon'}_{x,y} p_{j,y+\epsilon'}(\Psi^\star_{x})\,,\label{eq:cons22}\\
   &b_{x-1}(\Phi_j^\star) 
   p_{jy}(\Psi^\star_x)
=\frac{\nu(\Psi_x)}{\nu(\Psi_{x-1})} \sum_{\epsilon'=0,+1} D^{-,\epsilon'}_{x,y} p_{j,y+\epsilon'}(\Psi^\star_{x-1})\,.\label{eq:cons23}
\end{align}
\end{subequations}

\subsection{Resolution of the constraints  \eqref{eq:cons1} and \eqref{eq:cons2}\label{sec:cons}}

To solve the constraints  \eqref{eq:cons1} and \eqref{eq:cons2}, we pick 
all the possible choices for the parameter arrays $\Phi_j$ and $\Psi_x$ ($0\leq x,j \leq N$), recalled in Appendix \ref{app:pa}. Then, for each choice, we compute $c_{j+1}(\Psi_x)$, $a_j(\Psi_x)$, $b_{j-1}(\Psi_x)$. 
We want $a_j(\Psi_x)$ to be an affine transformation of the function depending on $x$ in the L.H.S. of 
the recurrence relation \eqref{eq:recu} for 
$p_{xi}(\Phi_j)$. Namely, we want 
\begin{equation}
    a_j(\Psi_x)= z^{0}_{j} t_x +u_j\,,
\end{equation}
where $z^0_j$ and $u_j$ are functions only of $j$ and we recall that $t_x$ are the eigenvalues associated to $\Phi_j$.
If $z^0_j$ and $u_j$ exist, we can compute $C^{\epsilon,0}_{i,j}$ and solve the constraints \eqref{eq:cons12}:
\begin{align}
   & C^{+,0}_{i,j}=z^{0}_{j}\ b_{i-1}(\Phi_j)\,,\\
   &  C^{-,0}_{i,j}=z^{0}_{j}\ a_i(\Phi_j)+u_j\,,\\
   &   C^{+,0}_{i,j}=z^{0}_{j}\ c_{i+1}(\Phi_j)\,.
\end{align}

For both other constraints \eqref{eq:cons11} and  \eqref{eq:cons13}, we demand that  
 $c_{j+1}(\Psi_x)$ and $b_{j-1}(\Psi_x)$ be the functions depending on $x$ in the L.H.S. of 
contiguity equations for $p_{xi}(\Phi_j)$. Such equations have the following forms:
\begin{subequations} \label{eq:contz}
    \begin{align}
   & \lambda_x^+(\Phi_j) p_{x,i}(\Phi_j)
    =\gamma^{0+}_i(\Phi_j) p_{x,i}(\Phi_{j+ 1})+\gamma^{-+}_i(\Phi_j)  p_{x,i-1}(\Phi_{j+1})\;,\\
   & \lambda_x^-(\Phi_j) p_{x,i}(\Phi_j)
    =\gamma^{+-}_i(\Phi_j) p_{x,i+1}(\Phi_{j-1})+\gamma^{0-}_i(\Phi_j)  p_{x,i}(\Phi_{j- 1})\;,
\end{align}
\end{subequations}
where $\lambda$ and $\gamma$ are functions to be determined.
Let us emphasize that in the R.H.S. of the previous equations there is a shift $j\pm 1$ for the parameter array in the argument of the polynomials.
In particular, this shift means that the polynomials on the R.H.S. and L.H.S. of these equations are associated to parameter arrays with different diameters. This type of equations is usually called contiguity relations. In Appendix \ref{app:pa}, such equations for the different parameter arrays are given. These equations are proved by using the method developed in \cite{icosi} for the Racah case. Below, we recall this method for the Krawtchouk case.
We want now that $c_{j+1}(\Psi_x)$ and $b_{j-1}(\Psi_x)$ correspond to 
$\lambda_x^+(j)$ and $\lambda_x^-(j)$ as follows 
\begin{align}
    c_{j+1}(\Psi_x) =z_j^{+}\lambda^+_x(\Phi_j) \,, \\
    b_{j-1}(\Psi_x)  =z_j^{-} \lambda^-_x(\Phi_j)\,, 
\end{align}
for some functions $z_j^{\pm}$ (independent of $x$).
Then the parameters $C^{\epsilon,\pm}_{i,j}$ can be computed
\begin{align}
    C^{\epsilon,\pm}_{i,j}=z_j^{\pm} \gamma_i^{\epsilon,\pm}(\Phi_j)\,,
\end{align} 
and the constraints \eqref{eq:cons1} are completely solved.

The same procedure is also done to solve the constraints \eqref{eq:cons2} and compute $D^{\epsilon,\epsilon'}_{x,y}$.\\

As just explained, it is important to obtain contiguity relations \eqref{eq:contz} for the considered polynomials. Here, we prove the contiguity recurrence equation  \eqref{eq:adkraw1} for the Krawtchouk polynomial. All the other contiguity equations given in Appendix \ref{app:pa} are proven similarly.

Relation \eqref{eq:adkraw1} is proven by induction. 
For $i=1$, it is checked by direct computation.
For a given $i\geq 1$, we suppose that the equation is valid for any $i'\leq i$.
From the recurrence relation \eqref{eq:recu}, one gets
\begin{equation}
p_{x,i+1}(\Phi)=\frac{1}{c_{i+1}(\Phi)}\Big((\theta_x-a_i(\Phi))p_{x,i}(\Phi)-b_{i-1}(\Phi)p_{x,i-1}(\Phi) \Big)\,.
\end{equation}
From the induction hypothesis, one may replace $p_{x,i}(\Phi)$ and $p_{x,i-1}(\Phi)$ to get 
\begin{align}
p_{x,i+1}(\Phi)=\frac{1}{c_{i+1}(\Phi)}\Big(&(\theta_x-a_i(\Phi))\big(p^-_{x,i}(\Phi) +\frac{r}{ss^\star-r}p^-_{x,i-1}(\Phi) \big)\nonumber\\
&-b_{i-1}(\Phi)\big(p^-_{x,i-1}(\Phi) +\frac{r}{ss^\star-r}p^-_{x,i-2}(\Phi) \big) \Big)\,,
\end{align}
where we recall that the superscript $-$ of $p^-_{x,i}(\Phi)$ means that $d$ is replaced by $d-1$ in its expression.
We use again twice the recurrence relation \eqref{eq:recu} 
for $\theta_x p^-_{x,i}(\Phi)$ and $\theta_x p^-_{x,i-1}(\Phi) $ to obtain
\begin{align}
p_{x,i+1}(\Phi)=\frac{1}{c_{i+1}(\Phi)}\Big(
&b_{i-1}^-(\Phi)p_{x,i-1}^-(\Phi)+a_{i}^-(\Phi)p_{x,i}^-(\Phi)+c_{i+1}^-(\Phi)p_{x,i+1}^-(\Phi)\nonumber\\
&+\frac{r}{ss^\star-r}\big( b_{i-2}^-(\Phi)p_{x,i-2}^-(\Phi)+a_{i-1}^-(\Phi)p_{x,i-1}^-(\Phi)+c_{i}^-(\Phi)p_{x,i}^-(\Phi)\big)\nonumber\\
&-a_i(\Phi)\big(p^-_{x,i}(\Phi) +\frac{r}{ss^\star-r}p^-_{x,i-1}(\Phi) \big)\nonumber\\
&-b_{i-1}(\Phi)\big(p^-_{x,i-1}(\Phi) +\frac{r}{ss^\star-r}p^-_{x,i-2}(\Phi) \big) \Big)\,,
\end{align}
where as for $p_{x,i}$, the superscript $-$ of $b$ and $c$ means that $d$ is replaced by $d-1$ in their expression. Then, using the explicit expressions \eqref{eq:abqkraw} for $a_i$, $b_i$ and $c_i$, one proves that relation \eqref{eq:adkraw1} is also satisfied for $i+1$ which concludes the proof of its validity for any $i$.

\subsection{Bilinear form}

If, for a choice of parameter arrays $\Phi_j$ and $\Psi_x$, the constraints \eqref{eq:cons1} and \eqref{eq:cons2} are solved, then
we want to find sufficient conditions on the parameter arrays for statement $(vii)$ of Definition \ref{def:IT} (concerning the bilinear form) to hold.

Let us choose a bilinear form such that 
\begin{equation}
    \langle V_{x,y},V_{u,v} \rangle=\delta_{x,u}\delta_{y,v} F(x,y)\,,
\end{equation}
for a non-vanishing function $F(x,y)$.

The change of basis \eqref{eq:VsV} allows us to compute
\begin{align}
    \langle V^\star_{i,j},V^\star_{k,\ell} \rangle&=
    \sum_{x=0}^{N-j}\sum_{y=0}^{N-x} p_{xi}(\Phi_j)p_{yj}(\Psi_x)
      \sum_{u=0}^{N-\ell}\sum_{v=0}^{N-u} p_{uk}(\Phi_\ell)p_{v\ell}(\Psi_u)
    \langle V_{x,y}, V_{u,v}\rangle \,,\\
    &= \sum_{x=0}^{N-j}\sum_{y=0}^{N-x} p_{xi}(\Phi_j)p_{yj}(\Psi_x)
     p_{xk}(\Phi_\ell)p_{y\ell}(\Psi_x)
   F(x,y)\,. \label{eq:sijkl}
\end{align}
The Wilson duality  (see for example \cite{TerClas}) states that 
\begin{align}\label{eq:Wil}
p_{ij}(\Phi)=\frac{k_j(\Phi)}{k_i(\Phi^\star)} p_{ji}(\Phi^\star)\,,\end{align}
where $k_i(\Phi)$ is given in \eqref{eq:kj}, which allows us to rewrite \eqref{eq:sijkl} as follows 
\begin{align}
    \langle V^\star_{i,j},V^\star_{k,\ell} \rangle&=
     \sum_{x=0}^{N-j}\sum_{y=0}^{N-x}
     \frac{k_i(\Phi_j)}{k_x(\Phi^\star_j)} 
      \frac{k_j(\Psi_x)}{k_y(\Psi^\star_x)} 
     p_{ix}(\Phi_j^\star)p_{jy}(\Psi^\star_x)
     p_{xk}(\Phi_\ell)p_{y\ell}(\Psi_x)
   F(x,y)\,. 
\end{align}
If there exists $G(i,j)$, a function which does not depend on $x$ and $y$, such that
\begin{align} \label{eq:FG}
     \frac{k_i(\Phi_j)}{k_x(\Phi^\star_j)} 
      \frac{k_j(\Psi_x)}{k_y(\Psi^\star_x)} F(x,y)=\frac{G(i,j)}{\nu(\Psi_x)\nu(\Phi_j)}\,, 
\end{align}
then, using two times the orthogonality relation \eqref{eq:Pinv}, one arrives at:
\begin{align}
    \langle V^\star_{i,j},V^\star_{k,\ell} \rangle&=\frac{G(i,j)}{\nu(\Phi_j)}
     \sum_{x=0}^{N-j}
           \frac{1}{\nu(\Psi_x)} 
     p_{ix}(\Phi_j^\star)  p_{xk}(\Phi_\ell) \sum_{y=0}^{N-x} p_{jy}(\Psi^\star_x)
   p_{y\ell}(\Psi_x)\\
   &=\delta_{i,k}\delta_{j,\ell} G(i,j)\,.
\end{align}
Finally, the statement $(vii)$ of Definition \ref{def:IT} of a factorized Leonard pair is satisfied if the constraint \eqref{eq:FG} is satisfied.

To conclude this section, let us summarize the results of this section in the following theorem.
\begin{theo}\label{thm:par}
The parameters arrays $\Phi_j$ and $\Psi_x$ given by 
\begin{subequations}
    \begin{align}
    &\Phi_j=(t_x ,t^\star_x(j),x=0,\dots,N-j;\ f_x(j), g_x(j),x=1,\dots,N-j)\,,\\
    &\Psi_x=(\theta_j(x) ,\theta^\star_j,j=0,\dots,N-x;\ \varphi_j(x), \phi_j(x),j=1,\dots,N-x)\,,
\end{align}
\end{subequations}
are called compatible if there exist $C_{i,j}^{\epsilon,\epsilon'}$ and $D_{x,y}^{\epsilon,\epsilon'}$ such that relations \eqref{eq:cons1} and \eqref{eq:cons2} are satisfied and the couple of functions $F(x,y)$ and $G(i,j)$ exists such that \eqref{eq:FG} is valid.

A factorized Leonard pair can be associated to a compatible pair of parameter arrays, and the coefficients $C_{i,j}^{\epsilon,\epsilon'}$ and $D_{x,y}^{\epsilon,\epsilon'}$ are uniquely determined. The actions of the elements of the basis $(X,Y;X^\star, Y^\star)$ of this factorized Leonard pair on the basis $V_{x,y}$ (statement (iii) of Definition \ref{def:IT}) are given by
\begin{subequations}
\begin{align}
&X V_{x,y}= t_x V_{x,y}\,,\qquad Y V_{x,y}= \theta_y(x) V_{x,y}\,,\\
&  X^\star V_{x,y} =\sum_{(\epsilon,\epsilon')\in A_2} D^{\epsilon,\epsilon'}_{x,y}V_{x+\epsilon,y+\epsilon'}\,,\\
&Y^\star V_{x,y}=b_{y-1}(\Psi_x^\star) V_{x,y-1}+
a_{y}(\Psi_x^\star) V_{x,y}+ c_{y+1}(\Psi_x^\star) V_{x,y+1}\,.
\end{align}
\end{subequations}
Similarly, on the basis $V^\star_{i,j}$ (statement (iv) of Definition \ref{def:IT}), the actions are given by 
\begin{subequations}
\begin{align}
&X V^\star_{i,j}=b_{i-1}(\Phi_j) V^\star_{i-1,j}+
 a_{i}(\Phi_j) V^\star_{i,j}+c_{i+1}(\Phi_j) V^\star_{i+1,j}\,,\\
 &  Y V^\star_{i,j} =\sum_{(\epsilon,\epsilon')\in A_2} C^{\epsilon,\epsilon'}_{i,j}V^\star_{i+\epsilon,j+\epsilon'}\,,\\
&X^\star V^\star_{i,j}= t^\star_i(j) V^\star_{i,j}\,,\qquad Y^\star V^\star_{i,j}= \theta_j^\star V^\star_{i,j}\,.
\end{align}
\end{subequations}
The change of basis between both bases is provided by \eqref{eq:VsV}.
\end{theo}

\section{Polynomials of Tratnik type \label{sec:pol}}

In this section, let $\Phi_j$ and $\Psi_x$ be compatible parameter arrays.
Before giving a list of possible factorized Leonard pairs, let us discuss the bivariate polynomials 
\begin{equation}\label{eq:Tij}
    T_{i,j}(x,y)= p_{xi}(\Phi_j)p_{yj}(\Psi_x) \,,
\end{equation}
which are the entries of the change of basis between $V_{x,y}$ and $V^\star_{i,j}$ given in \eqref{eq:VsV}.
These bivariate functions are called of Tratnik type in reference to \cite{Tra} where such polynomials have been introduced.

We shall show that these bivariate polynomials are bispectral, \textit{i.e.}\ they satisfy two recurrence and two difference relations, and that they are orthogonal.   

\subsection{Difference relations}

From the actions of $X^\star$ and $Y^\star$ on the vectors $V^\star_{i,j}$, one deduces the difference relations for $T_{i,j}(x,y)$. 
Indeed, using the results of Section \ref{sec:condtruction}, the action of $X^\star$ on $V^\star_{i,j}$ can be computed in two ways:
\begin{align}
    X^\star V_{i,j}^\star=t^\star_i(j)V_{i,j}^\star=\sum_{x=0}^{N-j}\sum_{y=0}^{N-x} T_{i,j}(x,y) t^\star_i(j) V_{x,y}\,,
\end{align}
and
\begin{align}
    X^\star V_{i,j}^\star&=\sum_{x=0}^{N-j}\sum_{y=0}^{N-x} T_{i,j}(x,y)  X^\star V_{x,y}\\
    &=\sum_{x=0}^{N-j}\sum_{y=0}^{N-x} T_{i,j}(x,y) \sum_{(\epsilon,\epsilon')\in A_2} D^{\epsilon,\epsilon'}_{x,y}V_{x+\epsilon,y+\epsilon'}\,.
\end{align}
Comparing the coefficients in front of $V_{x,y}$, one obtains the first difference relation
\begin{subequations}\label{eq:diffT}
    \begin{align}
   t^\star_i(j)  T_{i,j}(x,y)= \sum_{(\epsilon,\epsilon')\in A_2} D^{\epsilon,\epsilon'}_{x-\epsilon,y-\epsilon'}T_{i,j}(x-\epsilon,y-\epsilon')\,.
\end{align}
Similarly, the action of $Y^\star$ on $V^\star_{i,j}$ gives the second difference relation
\begin{align}
    \theta^\star_j T_{i,j}(x,y)=b_y(\Psi_x^\star)T_{i,j}(x,y+1)+
a_y(\Psi_x^\star)T_{i,j}(x,y)+
c_y(\Psi_x^\star)T_{i,j}(x,y-1)\,.
\end{align}
\end{subequations}

\subsection{Recurrence relations}

To obtain the recurrence relations, let us first invert the change of basis \eqref{eq:VsV}:
\begin{align}
V_{x,y}&=\sum_{j=0}^{N-x}\sum_{i=0}^{N-j}\frac{1}{\nu(\Phi_j)\nu(\Psi_x)} p_{ix}(\Phi^\star_j)  p_{jy}(\Psi_x^\star)V_{i,j}^\star\,.
\end{align}
Using the Wilson duality \eqref{eq:Wil}, one gets
\begin{align}
V_{x,y}&=\sum_{j=0}^{N-x}\sum_{i=0}^{N-j}\frac{   k_x(\Phi_j^\star)k_y(\Psi_x^\star)}{\nu(\Psi_x)\nu(\Phi_j)k_i(\Phi_j)k_j(\Psi_x)}  T_{i,j}(x,y) V_{i,j}^\star\,.
\end{align}
Relation \eqref{eq:FG} allows us to simplify the latter relation:
\begin{align}
V_{x,y}&=\sum_{j=0}^{N-x}\sum_{i=0}^{N-j}\frac{F(x,y)}{G(i,j) }  T_{i,j}(x,y) V_{i,j}^\star\,.\label{eq:VTVS}
\end{align}
One obtains recurrence relations for $T_{i,j}(x,y)$ from the actions of $X$ and $Y$ on the vectors $V_{x,y}$.
Indeed, the  action of $X$ on $V_{x,y}$ provides the first recurrence relation
\begin{subequations}\label{eq:recuT}
    \begin{align}\label{eq:recuT1}
    t_x T_{i,j}(x,y)=\frac{G(i,j)}{G(i+1,j)} b_i(\Phi_j) T_{i+1,j}(x,y)+
   a_i(\Phi_j) T_{i,j}(x,y)+
    \frac{G(i,j)}{G(i-1,j)} c_i(\Phi_j) T_{i-1,j}(x,y)\,,
\end{align}
and the  action of $Y$ on $V_{x,y}$ provides the second one
\begin{align}\label{eq:recuT2}
  \theta_y(x) T_{i,j}(x,y)=  \sum_{(\epsilon,\epsilon')\in A_2}\frac{G(i,j)}{G(i-\epsilon,j-\epsilon')} C^{\epsilon,\epsilon'}_{i-\epsilon,j-\epsilon'}T_{i-\epsilon,j-\epsilon'}(x,y)\,.
\end{align}
\end{subequations}

\subsection{Degree of the polynomials \label{sec:degree}}

Knowing that $p_{0x}(\Phi)=1$ for any parameter array $\Phi$ (see \cite{TerClas}), one deduces that $T_{0,0}(x,y)=1$. Then, using the recurrence relations \eqref{eq:recuT1} and \eqref{eq:recuT2} for $i,j=0$, one shows that $T_{1,0}(x,y)$ and $T_{0,1}(x,y)$ are polynomials with respect to $t_x$ and $\theta_y(x)$ of total degree $1$.

Suppose that, for a given $k$ $(0<k < N)$ and for $i+j<k$, $T_{i,j}(x,y)$ are bivariate polynomials with respect to $t_x$ and $\theta_y(x)$ of total degree $i+j$. Using the recurrence relation \eqref{eq:recuT1} for $(i,j)=(k-1,0),(k-2,1),\dots,(0,k-1)$, one can compute $T_{i,k-i}(x,y)$ for $i=1,2,\dots k$ in terms of $T_{i,j}(x,y)$ with $i+j< k$ and deduce that they are of total degree $k$. 
Then, using \eqref{eq:recuT2} for $(i,j)=(0,k-1)$, one deduces that $T_{0,k}(x,y)$ is also of total degree $k$ in terms of $t_x$ and $\theta_y(x)$. 
Therefore, one shows recursively with this method that $T_{i,j}(x,y)$ is a polynomial with respect to $t_x$ and $\theta_y(x)$ of total degree $i+j$, for any couple $(i,j)$.

In fact, it is possible to be a bit more precise. Following the same recursive proof as previously, one can also show that the degree of $T_{i,j}(x,y)$ with respect to $\theta_y(x)$ is exactly $j$.

\subsection{Orthogonality relation for the Tratnik polynomials\label{sec:ortho}}

The functions $T_{ij}(x,y)$ are also orthogonal in addition of being bispectral. Indeed, the vectors $V_{x,y}$ can be expressed in terms of the vectors $V^\star_{i,j}$ thanks to \eqref{eq:VTVS}, and the latter can in turn be replaced using \eqref{eq:VsV}. This leads to the expression
\begin{align}
V_{x,y}=\sum_{x'=0}^N \sum_{y'=0}^{N-x'} \ \sum_{j=0}^{\text{min}(N-x,N-x')}\sum_{i=0}^{N-j} \frac{F(x,y)}{G(i,j)} T_{i,j}(x,y)T_{i,j}(x',y')\ V_{x',y'}.
\end{align}
Identifying the coefficient in front of $V_{x,y}$, one gets the orthogonality relation
\begin{align}\label{eq:ortho}
\sum_{j=0}^{\text{min}(N-x,N-x')}\sum_{i=0}^{N-j} \frac{1}{G(i,j) } T_{i,j}
(x,y)T_{i,j}(x',y')= \delta_{x,x'}\delta_{y,y'}\ \frac{1}{F(x,y)}\, .
\end{align}

Now, starting by expressing $V^\star_{i,j}$ in terms of $V_{x,y}$ thanks to \eqref{eq:VsV}, and then 
replacing the latter using \eqref{eq:VTVS}, one gets
\begin{align}
\sum_{x=0}^{\text{min}(N-j,N-j')}\sum_{y=0}^{N-x} F(x,y) T_{i,j}
(x,y)T_{i',j'}(x,y)= \delta_{i,i'}\delta_{j,j'}\ G(i,j)\, .
\end{align}

\section{List of factorized $A_2$-Leonard pairs of type II \label{sec:listA2II}}

In the following sections, we provide couples of compatible parameter arrays, for $j=0,\dots, N$, 
\begin{subequations}
\begin{align}
&\Phi_j=(t_x,t^\star_x(j),x=0,\dots,N-j;\ f_x(j), g_x(j),x=1,\dots,N-j)\,,
\end{align}
and, for $x=0,\dots,N$,
\begin{align}
&\Psi_x=(\theta_j(x) ,\theta^\star_j,j=0,\dots,N-x;\varphi_j(x), \phi_j(x),j=1,\dots,N-x)\,.
\end{align}
\end{subequations}
Then, using Theorem \ref{thm:par}, this list also gives examples of factorized $A_2$-Leonard pairs.
We use the $A_2$-contiguity recurrence equations given in Appendix \ref{app:pa} to prove that the constraints \eqref{eq:cons1} and \eqref{eq:cons2} hold. We also verify that \eqref{eq:FG} is satisfied. 
We shall use the invariance \eqref{eq:inv} to simplify these parameter arrays. The procedure is detailed in the example of Section \ref{ssec:hahndhahn}.

We restrict in this section to the case where the parameter arrays are of type II. We recall that type II means that the (dual) eigenvalues can be written as $z_2 x^2+ z_1 x +z_0$, for some parameters $z_i$. Let us also recall that, for one type of eigenvalues, there exist different kinds of polynomials. One names a parameter array according to this associated polynomial.
The following conjecture provides a strong restriction on the choice of the parameter arrays.
\begin{conj}
To get a factorized Leonard pair, it is necessary that 
$\Phi_j$ for any $j$ (resp. $\Psi_x$ for any $x$) is  associated to the same kind of polynomial.
\end{conj}
\noindent
All the examples satisfy this conjecture and we name 
a pair of compatible parameter arrays $(\Phi_j,\Psi_x)$ with the two names 
of the associated polynomials.

We do not discuss the irreducibility of each Leonard pair. We refer to \cite{TerClas} where conditions on their parameters are explicitly given such that the associated Leonard pair is irreducible.

\subsection{Hahn / dual Hahn} \label{ssec:hahndhahn}

Let us consider the case where the parameter array $\Phi_j$ is of Hahn kind and $\Psi_x$ is of dual Hahn kind.
The parameter array $\Phi_j$ ($j=0,1,\dots,N$) reads as follows, for $x=0,1,\dots,N-j$,
\begin{subequations}\label{eq:paHdH1}
\begin{align}
&t_x= x,\qquad t^\star_x(j)= m^\star x(x+1+s^\star+j)\,,
\end{align}
and, for $x=1,\dots,N-j$,
\begin{align}
&f_x(j)= m^\star x(x-N+j-1)(x+r),\qquad
g_x(j)=- m^\star x(x-N+j-1)(x+s^\star+j-r)\,.
\end{align}
\end{subequations}
The parameter array $\Psi_x$ ($x=0,1,\dots,N$) reads as follows, for $j=0,1,\dots,N-x$,
\begin{subequations}\label{eq:paHdH2}
\begin{align}
&\theta_j(x)=  \mu j(j+1+\sigma+x),\qquad
\theta^\star_j=  j\,,
\end{align}
and, for $j=1,\dots,N-x$,
\begin{align}
& \varphi_j(x)= \mu j(j-N+x-1)(j+\rho),\qquad
\phi_j(x)=  \mu j(j-N+x-1)(j+\rho-\sigma-N-1)\,.
\end{align} 
\end{subequations}
The parameter array $\Phi_j$ is obtained from the one given in \eqref{eq:parhahn}, by taking
\begin{align}
\theta_0 \to 0, \quad \theta^\star_0 \to 0, \quad \sigma \to 1, \quad \mu^\star \to m^\star, \quad \sigma^\star \to s^\star +j, \quad \rho \to r, \quad d \to N-j \,,
\end{align}
and renaming $x\to i, \ \theta \to t, \varphi \to f, \ \phi \to g$. The dependence of $\Phi_j$ on $j$ is therefore present in the diameter $d$ and as a shift in the parameter $\sigma^\star$.
The parameters $\theta_0$ and $\theta^\star_0$ are sent to $0$ and $\sigma$ to $1$ without loss of generality due to the invariance \eqref{eq:inv}.
Similarly, the parameter array $\Psi_x$ is obtained from \eqref{eq:pardhahn} by taking
\begin{align}
\theta_0 \to 0, \quad \theta^\star_0 \to 0, \quad \sigma^\star\to 1, \quad \sigma \to  \sigma +x, \quad d \to N-x \,, 
\end{align}
and renaming $i \to j$. The dependence of $\Psi_x$ on $x$ is present in the diameter $d$ and as a shift in the parameter $\sigma$. 

The recurrence coefficients corresponding to the parameter array $\Phi_j$ are obtained from \eqref{eq:rechahn}:
\begin{subequations}
\begin{align}
a_i(\Phi_j)&= -b_{i}(\Phi_j) -c_{i}(\Phi_j) \,,\\
b_{i}(\Phi_j)&= \frac{(i+j-N)(i+1+r)(i+1+s^\star+j)}{(2i+j+s^\star+2)(2i+j+s^\star+1)}\,,\\
c_{i}(\Phi_j)&= -\frac{i(i+j+s^\star-r)(i+1+s^\star+N)}{(2i+j+s^\star)(2i+j+s^\star+1)}\,,
\end{align}
and those corresponding to the parameter array $\Psi_x$ are obtained from \eqref{eq:recdhahn}:
\begin{align}
a_j(\Psi_x)&= -b_j(\Psi_x) -c_j(\Psi_x) \,,\\
b_j(\Psi_x)&=\mu(j-N+x)(j+1+\rho)\,,\\
c_j(\Psi_x)&=\mu j(j-N-1+\rho-\sigma)\, .
\end{align}
\end{subequations}
The coefficient $a_j(\Psi_x)$ can be written as
\begin{equation}
a_j(\Psi_x) = -\mu(j+1+\rho)t_x +\mu(N-j)(j+1+\rho) + \mu j(N-j+1+\sigma-\rho) \,.
\end{equation}
The constraint \eqref{eq:cons12} can therefore be solved using the recurrence relation \eqref{eq:recu} of $p_{xi}(\Phi_j)$, giving
\begin{subequations}
\begin{align}
&C^{-,0}_{i,j}=-\mu(j+1+\rho) b_{i-1}(\Phi_j)\,, \\
&C^{+,0}_{i,j}=-\mu(j+1+\rho) c_{i+1}(\Phi_j)\,, \\
&C^{0,0}_{i,j}=-\mu(j+1+\rho) a_{i}(\Phi_j)+  \mu(N-j)(j+1+\rho) +\mu j(N-j+1+\sigma-\rho)\,.
\end{align}
The contiguity relation \eqref{eq:adhahn5} for $p_{xi}(\Phi_j)$ allows one to solve the constraint \eqref{eq:cons11}:
\begin{align}
&C^{-,+}_{i,j}=\frac{\mu(j+1)(j-N+\rho-\sigma)(i+r)(s^\star+j+2)}{(2i+s^\star+j)(s^\star+j+1-r)}\,, \\
&C^{0,+}_{i,j}=\frac{\mu(j+1)(j-N+\rho-\sigma)(i+1+s^\star+j-r)(s^\star+j+2)}{(2i+s^\star+j+2)(s^\star+j+1-r)}\,.
\end{align}
The constraint \eqref{eq:cons13} is solved using the contiguity relation \eqref{eq:adhahn6}: 
\begin{align}
&C^{0,-}_{i,j}=\frac{\mu(j+\rho)(N-j+1-i)(i+s^\star+j)(r-s^\star-j)}{(2i+s^\star+j)(s^\star+j+1)}\,, \\
&C^{+,-}_{i,j}=\frac{\mu(j+\rho)(i+1)(i+s^\star+N+2)(r-s^\star-j)}{(2i+s^\star+j+2)(s^\star+j+1)}\,.
\end{align}
\end{subequations}

In order to determine the coefficients $D^{\epsilon,\epsilon'}_{x,y}$ which appear in the constraints \eqref{eq:cons2}, one can observe using \eqref{eq:Phistar} that $\Phi_j^\star$ is the same as $\Psi_x$ with the replacements $j\to x, \ x \to j, \ \mu \to m^\star,\ \sigma \to s^\star, \ \rho \to r$, and $\Psi_x^\star$ is the same as $\Phi_j$ with the replacements $j\to x, \ x \to j, \ m^\star \to \mu,\ s^\star \to \sigma, \ r \to \rho$. Therefore, comparing the constraints \eqref{eq:cons1} and \eqref{eq:cons2}, one deduces that
\begin{equation}
    D^{\epsilon,\epsilon'}_{x,y} =  C^{\epsilon',\epsilon}_{y,x}\Big|_{(s^\star,r,\mu,\sigma,\rho)\to(\sigma,\rho,m^\star,s^\star,r)}\ \times\ \frac{\nu(\Psi_{x+\epsilon})}{\nu(\Psi_x)}.
\end{equation}

Using \eqref{eq:nu}, one finds
\begin{equation}
    \frac{\nu(\Psi_{x+1})}{\nu(\Psi_x)} = \frac{\sigma-\rho+x+1}{\sigma+x+2}\,, \qquad \frac{\nu(\Psi_{x-1})}{\nu(\Psi_x)} = \frac{\sigma+x+1}{\sigma-\rho+x}\,. 
\end{equation}
The constraint \eqref{eq:FG} is also satisfied and one gets
\begin{align}
   &F(x,y)=\binom{N}{x,y} 
\frac{(-1)^{N-y} (r+1)_x (\rho+1)_y (\rho-N-\sigma)_{N-x-y}}{(\sigma+x+y+1)_y(\sigma+x+2y+2)_{N-x-y}(r-s^\star-N)_x}\,,\\
    &G(i,j)=\binom{N}{i,j} 
\frac{ (r+1)_i[(2+s^\star+j)_{N-j}]^2(\rho+1)_j}{(s^\star+j-r+1,s^\star+j+1+i)_i(s^\star+j+2i+2)_{N-j-i}(\rho-\sigma-N)_j(r-s^\star-N)_{N-j} }\,,
\end{align}
where $\binom{N}{i,j}=\frac{N!}{(N-i-j)!i!j!}$ and $(x_1,\dots,x_j)_i$ corresponds to the Pochhammer symbols:
\begin{align}
  (x_1,\dots,x_j)_i =\prod_{k=0}^{i-1}(x_1+k)\dots  (x_j+k)\,.
\end{align}
Therefore, by Theorem \ref{thm:par}, we conclude that the pair $(\Phi_j,\Psi_x)$ with the parameter arrays \eqref{eq:paHdH1} and \eqref{eq:paHdH2} is associated to a factorized $A_2$-Leonard pair.

The bivariate polynomials \eqref{eq:Tij} associated to this factorized Leonard pair are obtained from the (univariate) Hahn polynomials \eqref{eq:polHahn} and dual Hahn polynomials \eqref{eq:poldHahn}. They read explicitly as follows, for $i,j,x,y$ non-negative integers, $i+j\leq N$, $x+y\leq N$, 
\begin{align}\label{eq:TratHahn}
    T_{i,j}(x,y)= &\binom{N-j}{i} \binom{N-x}{j}
\frac{(-1)^j (r+1)_i(2+s^\star+j)_{N-j}(\rho+1)_j}{(s^\star+j-r+1,s^\star+j+1+i)_i(s^\star+j+2i+2)_{N-j-i}(\rho-\sigma-N)_j} \nonumber \\
& \times \  {}_3F_2 \Biggl({{-i,\;i+1+s^\star+j, \;-x}\atop
{r+1,\; -N+j}}\;\Bigg\vert \; 1\Biggr) \ {}_3F_2 \Biggl({{-j, \;-y,\; y+1+\sigma+x}\atop
{\rho+1,\; -N+x}}\;\Bigg\vert \; 1\Biggr)\,.
\end{align}
In the previous equation, ${}_3F_2$ corresponds to the usual hypergeometric function (see \textit{e.g.}\ \cite{Koek}).
As explained in Section \ref{sec:degree}, $T_{i,j}(x,y)$, given in the previous relation, are bivariate polynomials of variables 
$x$ and $\mu y(y+1+\sigma+x)$ of total degree $i+j$.
Let us emphasize that the functions obtained here are not the ones studied previously for example in \cite{Tra, GI}. 
The latter are not obtained here since one of the four bispectral relations are more complicated and involve nine terms (called $B_2$ in the outlooks of the article).

\subsection{Krawtchouk / dual Hahn}

Let us consider the case Krawtchouk / dual Hahn where $\Phi_j$ is given by, for $j=0,1,\dots,N$,
\begin{subequations}\label{eq:paKdH1}
\begin{align}
&t_x= x,\ 
t_x^\star(j)= x\,,\\ 
&f_x(j)=r x(x-N+j-1),\  
g_x(j)=(r -1) x(x-N+j-1)\,,
\end{align}
\end{subequations}
and, $\Psi_x$ is given by \eqref{eq:paHdH2}.

The coefficients $C^{\epsilon,\epsilon'}_{i,j}$ can be computed from \eqref{eq:adkraw1}-\eqref{eq:adkraw2} and $D^{\epsilon,\epsilon'}_{x,y}$, from \eqref{eq:adhahn5}-\eqref{eq:adhahn6}. Alternatively, one can take the results of Subsection \ref{ssec:hahndhahn}, perform the replacements $s^\star \to \frac{1}{m^\star}, \ r \to \frac{r}{m^\star}$, and then take the limit $m^\star \to 0$.

\subsection{Hahn / Krawtchouk}

Let us consider the case Hahn / Krawtchouk where $\Phi_j$ is given by  \eqref{eq:paHdH1} 
and $\Psi_x$ is given by, for $x=0,1,\dots,N$,
\begin{subequations}\label{eq:paHK2}
\begin{align}
&\theta_j(x)=  j,\ 
\theta^\star_j= j\,,\\ 
&\varphi_j(x)=\rho j(j-N+x-1),
\phi_j(x)= (\rho-1)j(j-N+x-1)\,.
\end{align} 
\end{subequations}

The coefficients $C^{\epsilon,\epsilon'}_{i,j}$ can be computed as previously from \eqref{eq:adhahn5}-\eqref{eq:adhahn6} and $D^{\epsilon,\epsilon'}_{x,y}$, from \eqref{eq:adkraw1}-\eqref{eq:adkraw2}. Alternatively, one can take the results of Subsection \ref{ssec:hahndhahn}, perform the replacements $\sigma \to \frac{1}{\mu}, \ \rho \to \frac{\rho}{\mu}$, and then take the limit $\mu \to 0$.

The functions of  Tratnik type for the last two subsections (obtained as limit of the ones given in \eqref{eq:TratHahn}) have been used in \cite{TAG,BCPVZ,CVZZ} to study the non-binary Johnson association scheme.

\subsection{Krawtchouk / Krawtchouk \label{sec:kk}}
Let us consider the case Krawtchouk / Krawtchouk where $\Phi_j$ is given by \eqref{eq:paKdH1} and $\Psi_x$ is given by \eqref{eq:paHK2}. 

The coefficients $C^{\epsilon,\epsilon'}_{i,j}$ and $D^{\epsilon,\epsilon'}_{x,y}$ can be computed from \eqref{eq:adkraw1}-\eqref{eq:adkraw2}. Alternatively, one can take the results of Subsection \ref{ssec:hahndhahn}, perform the replacements $s^\star \to \frac{1}{m^\star}, \ r \to \frac{r}{m^\star}, \ \sigma \to \frac{1}{\mu}, \ \rho \to \frac{\rho}{\mu}$, and then take the limits $m^\star \to 0, \ \mu \to 0$.
 
The polynomial of Tratnik type in this case can be obtained as a limit of \eqref{eq:TratHahn}. This polynomial has already appeared in different contexts and it is a particular case of the Griffiths (or sometimes called Rahman) polynomials \cite{Ros,Zhe97,Gru,HR,MT,IT,GVZ,CVV,CFG}. 
Indeed, up to a normalisation and specifying some parameters, the Griffiths polynomials can be written as a product of two Krawtchouk polynomials and take the form \eqref{eq:Tij}.

\section{List of factorized $A_2$-Leonard pairs of type I\label{sec:listA2I}}

In the following, we provide couples of parameter arrays $(\Phi_j,\Psi_x)$ of type I, such that the constraints \eqref{eq:cons1}, \eqref{eq:cons2} and \eqref{eq:FG} can be solved. We recall that type I means that the (dual) eigenvalues can be written as $z^+ q^x+ z^- q^{-x} +z^0$, for some parameters $z^\pm,z^0$. These provide examples of factorized $A_2$-Leonard pairs of type I.

\subsection{ $q$-Hahn  /  dual $q$-Hahn } \label{ssec:qHdqH}

One can consider for $\Phi_j$ the parameter array of the $q$-Hahn kind
\begin{subequations}\label{ex:qhahn}
    \begin{align} 
&t_x= q^{-x}-1\,,\ t^\star_x(j)=m^\star q^{N-j}(q^{-x}-1)(1-s^\star q^{j+x+1}) \,,\\
&f_x(j)= m^\star q^{N-j+1-2x}(1-q^x)(1-q^{x-N+j-1})(1-r q^x) \,,\\
&g_x(j)=- m^\star q^{N-j+1-x}(1-q^x)(1-q^{x-N+j-1})(r-s^\star q^{x+j}) \,,
\end{align}
\end{subequations}
and for $\Psi_x$ one takes the parameter array of the dual $q$-Hahn kind
\begin{subequations}\label{ex:dqhahn}
    \begin{align}
&\theta_j(x)= \mu q^{N-x}(q^{-j}-1)(1-\sigma  q^{x+j+1})\,,
\ \theta^\star_j=q^{-j}-1 \,,\\
&\varphi_j(x)= \mu q^{N-x+1-2j}(1-q^j)(1-q^{j-N+x-1})(1-\rho q^j) \,,\\
& \phi_j(x)= \mu q^{2N-2x+2-2j}(1-q^j)(1-q^{j-N+x-1})(\sigma q^x -\rho q^{j-N+x-1}) \,.
\end{align}
\end{subequations}
The recurrence coefficients for these parameter arrays are obtained from \eqref{eq:recqH} and \eqref{eq:recdqH}:
\begin{subequations}
\begin{align}
a_i(\Phi_j)&= -b_{i}(\Phi_j)-c_{i}(\Phi_j) \,,\\
b_{i}(\Phi_j)&= \frac{(1-s^\star q^{i+j+1})(1-q^{i-N+j})(1-rq^{i+1})}{(1-s^\star q^{2i+j+2})(1-s^\star q^{2i+j+1})}\,,\\
c_{i}(\Phi_j)&= \frac{q^{i-N+j}(1-s^\star q^{N+i+1})(1-q^i)(s^\star q^{i+j}-r)}{(1-s^\star q^{2i+j+1})(1-s^\star q^{2i+j})}\,,
\end{align}
and
\begin{align}
a_j(\Psi_x)&= -b_j(\Psi_x) - c_j(\Psi_x) \,,\\
b_j(\Psi_x)&=\mu q^{N-x}(1-q^{j-N+x})(1-\rho q^{j+1})\,,\\
c_j(\Psi_x)&=\mu q^{N+1}(1-q^j)(\sigma-\rho q^{j-N-1})\, .
\end{align}
One can write
\begin{equation}
    a_j(\Psi_x)= -\mu q^{N}(1-\rho q^{j+1})t_x+\mu q^j(1-q^{N-j})(1-\rho q^{j+1})+\mu q^j(1-q^j)(\rho-\sigma q^{N-j+1}) \,.
\end{equation}
\end{subequations}
Using these recurrence coefficients, the recurrence relation \eqref{eq:recu} and the contiguity recurrence relations \eqref{eq:adqh1}--\eqref{eq:adqh2} for the polynomials $p_{xi}(\Phi_j)$, one solves the constraints \eqref{eq:cons1}:
\begin{subequations}
\begin{align}
&C^{-,0}_{i,j}=-\mu q^{N}(1-\rho q^{j+1}) b_{i-1}(\Phi_j)\,, \\
&C^{0,0}_{i,j}=-\mu q^{N}(1-\rho q^{j+1}) a_{i}(\Phi_j)+  \mu q^j(1-q^{N-j})(1-\rho q^{j+1}) +\mu q^j(1-q^j)(\rho-\sigma q^{N-j+1})\,, \\
&C^{+,0}_{i,j}=-\mu q^{N}(1-\rho q^{j+1}) c_{i+1}(\Phi_j)\,, \\
&C^{-,+}_{i,j}=\frac{\mu q^{N}(1-q^{j+1})(\sigma-\rho q^{j-N})(1-rq^i)(1-s^\star q^{j+2})}{(1-s^\star q^{j+2i})(r-s^\star q^{j+1})}\,, \\
&C^{0,+}_{i,j}=\frac{\mu q^{N+i+1}(1-q^{j+1})(\sigma-\rho q^{j-N})(r-s^\star q^{j+i+1})(1-s^\star q^{j+2})}{(1-s^\star q^{j+2i+2})(r-s^\star q^{j+1})}\,, \\
&C^{0,-}_{i,j}=-\frac{\mu q^{i+j}(1-\rho q^j)(1-q^{N-j+1-i})(1-s^\star q^{i+j})(r-s^\star q^{j})}{(1-s^\star q^{2i+j})(1-s^\star q^{j+1})}\,, \\
&C^{+,-}_{i,j}=-\frac{\mu q^j(1-\rho q^j)(1-q^{i+1})(1-s^\star q^{i+N+2})(r-s^\star q^j)}{(1-s^\star q^{2i+j+2})(1-s^\star q^{j+1})}\,.
\end{align}
\end{subequations}
For similar reasons as in Subsection \ref{ssec:hahndhahn}, the coefficients solving the constraints \eqref{eq:cons2} are given by
\begin{equation}
    D^{\epsilon,\epsilon'}_{x,y} =  C^{\epsilon',\epsilon}_{y,x}\Big|_{(s^\star,r,\mu,\sigma,\rho)\to(\sigma,\rho,m^\star,s^\star,r)}\ \times\ \frac{\nu(\Psi_{x+\epsilon})}{\nu(\Psi_x)}\, ,
\end{equation}
where
\begin{equation}
    \frac{\nu(\Psi_{x+1})}{\nu(\Psi_x)} = \frac{\sigma q^x-\rho q^{-1}}{\sigma q^x-q^{-2}}\,, \qquad \frac{\nu(\Psi_{x-1})}{\nu(\Psi_x)} = \frac{\sigma q^{x}-q^{-1}}{\sigma q^{x}-\rho}\,. 
\end{equation}
The constraint \eqref{eq:FG} is satisfied with the following functions 
\begin{align}
    F(x,y) = 
    \left[\begin{array}{c}N\\x,y\end{array}\right]_q
    \frac{(-1)^{N} q^{y(x+y-N)}(r q;q)_x(\rho q;q)_y(\rho q^{-N}/\sigma;q)_{N-x-y}}
    {(qr)^x(q^{N-x+1}s^\star/r;q)_x(\sigma q^{x+y+1};q)_y( q^{-y-N-1}/\sigma;q)_{N-x-y}} \,, 
\end{align}
\begin{align}
    G(i,j) = \left[\begin{array}{c}N\\i,j\end{array}\right]_q
    \frac{q^{-\frac{N}{2}(N+3)}q^{j(j+2)}(rq;q)_i(\rho q;q)_j}{(rq)^i(\sigma q^{N+2})^j(s^\star q^{j+1}/r,s^\star q^{i+j+1};q)_i(\rho q^{-N}/\sigma;q)_j} \nonumber \\
    \times \frac{(s^\star q^{j+2};q)^2_{N-j}}{(s^\star)^{N-j}(r q^{-N}/s^\star;q)_{N-j}(s^\star q^{2i+j+2};q)_{N-i-j}}\,,
\end{align}
where we are using standard notations for the $q$-Pochhammer symbols and the $q$-trinomial coefficient:
\begin{gather}
(x;q)_i=\prod_{k=0}^{i-1} (1-xq^k), \quad (x_1,\dots,x_j;q)_i = \prod_{k=1}^{j} (x_k;q)_i\,,\\
\left[\begin{array}{c}N\\i,j\end{array}\right]_q=\frac{(q;q)_N}{(q;q)_i(q;q)_j(q;q)_{N-i-j}}\,.
\end{gather}
By Theorem \ref{thm:par}, we conclude that the pair $(\Phi_j,\Psi_x)$ with the parameter arrays \eqref{ex:qhahn} and \eqref{ex:dqhahn} is associated to a factorized $A_2$-Leonard pair.

The bivariate polynomials \eqref{eq:Tij} associated to this factorized Leonard pair are obtained from the (univariate) $q$-Hahn polynomials \eqref{eq:polqHahn} and dual $q$-Hahn polynomials \eqref{eq:poldqHahn}, for $i,j,x,y$ non negative integers, $i+j\leq N$, $x+y\leq N$:
\begin{align}\label{eq:dqHH}
    &T_{i,j}(x,y)= \frac{(-1)^j q^{\frac{j}{2}(j-3-2N)}(rq;q)_i (s^\star q^{j+2} ;q)_{N-j}(\rho q;q)_j}
{(rq)^i(s^\star q^{j+1}/r,s^\star q^{i+j+1};q)_i (s^\star q^{2+2i+j};q)_{N-i-j}\sigma^j (q^{-N} \rho/\sigma;q)_j  } \nonumber \\
&\times \left[\begin{array}{c}N-j\\i\end{array}\right]_q \left[\begin{array}{c}N-x\\j\end{array}\right]_q \ {}_3\phi_2 \Biggl({{q^{-i},\; q^{-x}, \; s^\star q^{i+j+1} }\atop
{rq, \; q^{-N+j}}}\;\Bigg\vert \; q,q\Biggr) \ {}_3\phi_2 \Biggl({{q^{-j}, \; q^{-y},\;\sigma q^{x+y+1}}\atop
{\rho q, \; q^{-N+x}}}\;\Bigg\vert \; q,q\Biggr)\,.
\end{align}
As explained in Section \ref{sec:degree}, $T_{i,j}(x,y)$, given in the previous relation, are bivariate polynomials of variables 
$q^{-x}-1$ and $ \mu q^{N-x}(q^{-y}-1)(1-\sigma  q^{x+y+1})$ of total degree $i+j$.

\subsection{Limits}
As for the factorized Leonard pairs of type II in Section \ref{sec:listA2II}, other examples of factorized Leonard pairs of type I can be obtained by taking limits of the case $q$-Hahn/dual $q$-Hahn described in Subsection \ref{ssec:qHdqH}. We provide here only the limiting procedures leading to couples $(\Phi_j,\Psi_x)$ of different kinds; it is then straightforward to deduce the explicit expressions for the parameter arrays $\Phi_j,\Psi_x$, the coefficients $C_{i,j}^{\epsilon,\epsilon'},D_{x,y}^{\epsilon,\epsilon'}$, the functions $F(x,y),G(i,j)$ and the polynomials $T_{i,j}(x,y)$ using the results of Subsection \ref{ssec:qHdqH}.

For the parameter array $\Phi_j$, one can take the following limits starting from the $q$-Hahn case given by \eqref{ex:qhahn}:
\begin{itemize}[itemsep=0pt]
    \item quantum $q$-Krawtchouk: $s^\star \to \frac{q^{-N-1}}{m^\star}, \ r \to  \frac{1}{r m^\star}$ and then $m^\star \to 0$;
    \item affine $q$-Krawtchouk: $m^\star \to 1, \  s^\star \to 0$;
    \item $q$-Krawtchouk: $r \to 0$.
\end{itemize}
Similarly, for the parameter array $\Psi_x$, one can take the following limits starting from the dual $q$-Hahn case given by \eqref{ex:dqhahn}:
\begin{itemize}[itemsep=0pt]
    \item dual quantum $q$-Krawtchouk: $\sigma \to \frac{q^{-N-1}}{\mu}, \ \rho \to \frac{1}{\rho \mu}$ and then $\mu \to 0$;
    \item affine $q$-Krawtchouk: $\mu \to 1, \  \sigma \to 0$;
    \item dual $q$-Krawtchouk: $  \rho \to 0$.
\end{itemize}
Any combination of the above limits for $\Phi_j$ and $\Psi_x$ (including the case where no limit is taken for one or both of the parameter arrays) gives a factorized $A_2$-Leonard pair of type I.

The polynomial associated to the case quantum $q$-Krawtchouk / dual quantum $q$-Krawtchouk has been studied previously \cite{Tra,Rosengren, Sca, GPV}. In these previous papers (see \textit{e.g.}\ \cite{GPV}), it has been written as proportional to, for $n_1,\ n_2,\ x_1,\ x_2$ non negative integers with $n_1+n_2\leq N$ and $x_1+x_2\leq N$:
\begin{align}
  \mathbf{K}_{n_1n_2}(x_1,x_2)\propto 
   {}_2\phi_1 \Biggl({{\tilde{q}^{-n_1}, \;\tilde{q}^{-x_1}}\atop
{\tilde{q}^{-x_1-x_2} }}\;\Bigg\vert \; \tilde{q},\frac{\tilde{q}^{n_1+1}}{\alpha_1^2}\Biggr)
{}_2\phi_1 \Biggl({{\tilde{q}^{-n_2}, \;\tilde{q}^{-x_1-x_2+n_1}}\atop
{\tilde{q}^{-N+n_1} }}\;\Bigg\vert \; \tilde{q},\frac{\tilde{q}^{n_2+1}}{\alpha_2^2}\Biggr)\,. \label{eq:KT}
\end{align}
To recover the expression obtained in this paper, one performs the following changes:
\begin{align}
    x_1=y\,, \quad x_2=N-y-x\,, \quad n_1=j\quad\ n_2=N-i-j\,,\quad
    \alpha_1^2=1/r\,,\quad \alpha_2^2=1/\rho\,,
\end{align}
and, using relation (1.4.6) in \cite{GR}, one rewrites the second factor to get
\begin{align}
 \mathbf{K}_{j,N-i-j}(y,N-x-y)\propto    {}_2\phi_1 \Biggl({{\tilde{q}^{-j}, \;\tilde{q}^{-y}}\atop
{\tilde{q}^{-N+x} }}\;\Bigg\vert \; \tilde{q},r\tilde{q}^{j+1}\Biggr)\ 
{}_2\phi_1 \Biggl({{\tilde{q}^{-i}, \;\tilde{q}^{-x}}\atop
{\tilde{q}^{-N+j} }}\;\Bigg\vert \; \tilde{q},\rho \tilde{q}^{x+1}\Biggr)\,.
\end{align}
Then, setting $\tilde q=q^{-1}$ and using 
\begin{align}
     {}_2\phi_1 \Biggl({{a_1, \;a_2}\atop
{b }}\;\Bigg\vert \; q^{-1},z\Biggr)= {}_2\phi_1 \Biggl({{1/a_1, \;1/a_2}\atop
{1/b }}\;\Bigg\vert \; q,\frac{za_1a_2}{qb}\Biggr) \,,
\end{align}
one obtains 
\begin{align}\label{eq:KK}
 \mathbf{K}_{j,N-i-j}(y,N-x-y)\propto   {}_2\phi_1 \Biggl({{q^{-j}, \;{q}^{-y}}\atop
{{q}^{-N+x} }}\;\Bigg\vert \; q,r{q}^{j+x-N}\Biggr)\ 
{}_2\phi_1 \Biggl({{{q}^{-i}, \;{q}^{-x}}\atop
{{q}^{-N+j} }}\;\Bigg\vert \; {q},\rho{q}^{x+j-N}\Biggr)\,.
\end{align}
In \ref{eq:KK}, we recognize our expression of \eqref{eq:dqHH} in the limit quantum $q$-Krawtchouk / dual quantum $q$-Krawtchouk.

\section{Factorized $A_M$-Leonard pair }\label{sec:introM}

Up to now, we focused on the $A_2$-Leonard pairs. In this section, we propose the definition of $A_M$-Leonard pair for $M>2$ which is, as for the case $M=2$, a refinement of the definition of the rank $M$ Leonard pair (called here $A_M$-Leonard pair) given in \cite{IT}.

Let $N,M$ be positive integers and let $\cD$ be the set constituted of $M$-tuples of non-negative integers whose sum is at most $N$:
\begin{equation}
 \cD=\{ (x_1,x_2,\dots, x_M)\ | \ x_1,x_2,\dots, x_M \in \bZ_{\geq 0}, \ x_1+x_2+\dots+ x_M \leq N \}\,.
\end{equation}
Two $M$-tuples $\bold{x}=(x_1,x_2,\dots, x_M),\ \bold{y}=(y_1,y_2,\dots, y_M) \in \cD$ are called $A_M$-adjacent, $\bold{x}\ \adj\ \bold{y}$, if
\begin{equation}
  (x_1-y_1,x_2-y_2,\dots, x_M-y_M) \in A_M,
\end{equation}
where $A_M$ is the following subset of $\bZ^M$:
\begin{align}
 &A_M=\{(0,0,\dots,0),(1,-1,0,\dots,0), (1,0,\dots,0),  (-1,0,\dots,0),\text{ permutations } \}\,.
 \end{align}
The notation $A_M$ for the previous set comes from the fact that it corresponds to the root system of the Weyl group $A_M$.

\begin{defi}\cite{IT}
Let $V$ be a $\mathbb{C}$-vector space.
The pair $(H,H^\star)$ is a $A_M$-Leonard pair on the domain $\cD$ if the following statements are satisfied:
\begin{itemize}
 \item[(i)] $H$ is a $M$-dimensional  subspace of $\text{End}(V)$ whose elements are diagonalizable and mutually commute;
 \item[(ii)] $H^\star$ is a $M$-dimensional  subspace of $\text{End}(V)$ whose elements are diagonalizable and mutually commute;
 \item[(iii)] there exists a bijection $\bold{x}\mapsto V^{(M)}_\bold{x}$  from $\cD$ to the common eigenspaces of $H$ such that, for all $\bold{x} \in \cD$,
 \begin{equation}\label{eq:HsVM}
  H^\star V^{(M)}_\bold{x} \subseteq \sum_{ \genfrac{}{}{0pt}{}{\bold{y} \in \cD}{\bold{y}\ \adj \  \bold{x}}} V^{(M)}_\bold{y}\,;
 \end{equation}
 \item[(iv)]  there exists a bijection $\bold{x}\mapsto V^{(0)}_\bold{x}$  from $\cD$ to the common eigenspaces $V^{(0)}_\bold{x}$ of $H^\star$ such that, for all $\bold{x}\in \cD$,
 \begin{equation}\label{eq:HVsM}
  H V^{(0)}_\bold{x} \subseteq \sum_{ \genfrac{}{}{0pt}{}{\bold{y} \in \cD}{\bold{y}\ \adj \  \bold{x} }} V^{(0)}_\bold{y}\,;
 \end{equation}
 \item[(v)] there does not exist a subspace $W$ of $V$ such that $HW\subseteq W$, $H^\star W\subseteq W$, $W\neq 0$, $W\neq V$;
 \item[(vi)] $\text{dim}(V^{(0)}_\bold{x})=\text{dim}(V^{(M)}_\bold{x})=1$, for all  $\bold{x}\in \cD$;
 \item[(vii)] there exists a non-degenerate symmetric bilinear form  $\langle,\rangle$ on $V$ such that, for $\bold{x},\bold{y}\in \cD$ and $\bold{x} \neq \bold{y}$,
 \begin{equation}
\langle V^{(0)}_\bold{x}, V^{(0)}_\bold{y}\rangle=0\ ,\qquad  \langle V^{(M)}_\bold{x}, V^{(M)}_\bold{y}\rangle=0\,.
 \end{equation}
\end{itemize}
\end{defi}
Let us notice that the dimension of the vector space $V$ is finite and is given more precisely by $\dim(V) = |\cD|$. 

As previously for the case $A_2$, we want to add some statements to this previous definition in order to obtain a definition that is easier to handle.
\begin{defi} \label{def:flpM}
Let $(H,H^\star)$ be a $A_M$-Leonard pair on the domain $\cD$ and $X_1,X_2,\dots,X_M\in \text{End}(V)$ (resp. $X^\star_1,X^\star_2,\dots,X^\star_M$) be a basis of $H$ (resp. $H^\star$).
The pair $(H,H^\star)$ is a \textit{factorized $A_M$-Leonard pair} if the following statements are satisfied:
\begin{itemize}
 \item[(viii)]For $a=1,\dots,M-1$ and $\bold{x}=(x_1,\dots,x_M)$, the eigenvalues of $X_{a}$ on $V^{(M)}_\bold{x}$ 
 do not depend on $x_{a+1},x_{a+2},\dots x_M$:
 \begin{align}
 &\Big( X_{a}- \theta^{(a)}(x_1,\dots,x_a) \un\Big)V^{(M)}_\bold{x}=0\,,
 \end{align}
 and $\theta^{(a)}(x_1,\dots,x_{a-1},x_a)\neq \theta^{(a)}(x_1,\dots,x_{a-1},x'_a)$ for $x_a\neq x'_a$.\\
For $a=2,\dots,M$  and $\bold{i}=(i_1,\dots,i_M)$, the eigenvalues of $X^\star_{a}$ on $V^{(0)}_\bold{i}$ 
 do not depend on $i_{1},i_{2},\dots i_{a-1}$:
 \begin{align}
 &\Big( X^\star_{a}- {\theta^\star}^{(a)}(i_a,\dots,i_M)\un\Big) V^{(0)}_\bold{i}=0\,,
 \end{align}
 and ${\theta^\star}^{(a)}(i_a,i_{a+1},\dots,i_M)\neq {\theta^\star}^{(a)}(i'_a,i_{a+1},\dots,i_M)$ for $i_a\neq i'_a$;
  \item[(ix)] for $1\leq i < j \leq M$, one has
 \begin{align}\label{eq:com12}
  [X_i,X_j^\star]=0\,;
 \end{align}
 \item[(x)] let $1\leq a \leq M-1$. There exists a bijection $\bold{x}\mapsto V^{(a)}_\bold{x}$  from $\cD$ to the common eigenspaces of $X_1,X_2,\dots, X_{a}$ and $X_{a+1}^\star,\dots, X^\star_M$.
 \item[(xi)] For a given $b \in \{1,2\dots,M-1$\} and for any $a$ and $c$ satisfying $0< a \leq b <c \leq M$, the following relations hold, for $(\dots, i_a,\dots), (\dots, i_a+1,\dots)\in \cD$,
 \begin{subequations}
  \begin{align}
 & \langle V^{(a-1)}_{(\dots,i_a+1,\dots)} ,X_a V^{(a-1)}_{(\dots, i_a,\dots)}\rangle \neq 0\,,&& 
 \langle V^{(a-1)}_{(\dots,i_a,\dots)} , X_a V^{(a-1)}_{(\dots, i_a+1,\dots)}\rangle \neq 0\,, \label{eq:irre12}
 \end{align}
 and, for $(\dots, i_c,\dots), (\dots, i_c+1,\dots)\in \cD$,
 \begin{align}
      &\langle V^{(c)}_{(\dots,i_c+1,\dots)}  , X_c^\star V^{(c)}_{(\dots,i_c,\dots)}\rangle \neq 0\,,&& \langle V^{(c)}_{(\dots,i_c,\dots)} ,X_c^\star V^{(c)}_{(\dots,i_c+1,\dots)}\rangle \neq 0\,.\label{eq:irre22}
  \end{align}
  \end{subequations}
 \end{itemize}
 The ordered tuple $(X_1,X_2,\dots, X_M;X^\star_1,X^\star_2,\dots, X^\star_M)$ is called a basis of the factorized $A_M$-Leonard pair.
 \end{defi}

We believe that for a generic $M$ it is possible to obtain the same results as those obtained in this paper for the case $M=2$. In particular the 
polynomials of Tratnik type with $M$ variables may appear naturally. We plan to study the case with a generic $M$ in the future.

\section{Outlooks}\label{sec:conclu}

In this paper, we provide the definition of a factorized $A_2$-Leonard pair and start the study of these objects. We show that, in comparison to generic rank 2 Leonard pairs, it is possible to derive interesting properties for factorized Leonard pairs as a consequence of their refined definition, and it is possible to find several classes of examples based on combinations of polynomials of the ($q$-)Askey scheme. 
Regarding the examples, let us mention the following conjecture.
\begin{conj}
    The list of factorized $A_2$-Leonard pairs given in Sections \ref{sec:listA2II} and \ref{sec:listA2I} is exhaustive.
\end{conj}
\noindent In other words, we believe that a classification of all the factorized $A_2$-Leonard pairs is obtained. In particular, it seems that there is no factorized $A_2$-Leonard pairs involving the Bannai--Ito polynomials (of type III) or the ($q$-)Racah polynomials.      

To prove the above conjecture, we intend to study in detail the algebra associated to a factorized $A_2$-Leonard pair, mimicking the proofs for the usual Leonard pairs. Indeed, as shown in \cite{TV03}, the two members of a Leonard pair satisfy the Askey--Wilson relations. Using this result, we can show that, 
for a basis $(X,Y;X^\star,Y^\star)$ of a factorized $A_2$-Leonard pair, 
$X$ and $X^\star$ also satisfy these relations \textit{i.e.}
\begin{subequations}
\begin{align}
   & X^2 X^\star -\beta X X^\star X  +X^\star X^2 -\gamma (X X^\star + X^\star X) -\rho X^\star= \gamma^\star X^2+\omega X+ \eta I\,, \\
  &    (X^\star)^2 X -\beta X^\star X X^\star  +X (X^\star)^2 -\gamma^\star (X X^\star + X^\star X) -\rho^\star X= \gamma (X^\star)^2+\omega X^\star+ \eta^\star I\,,
\end{align}
\end{subequations}
where the parameters $\beta$, $\gamma$, $\gamma^\star$, $\rho$, $\rho^\star$, $\omega$, $\eta$ and $\eta^\star$ depend on $Y^\star$. We recall that $Y^\star$ commutes with $X$ and $X^\star$, meaning that $Y^\star$ is a central element in the algebra generated by $X$ and $X^\star$.
Similarly, $Y$ and $Y^\star$ also satisfy the Askey--Wilson relations:
\begin{subequations}
\begin{align}
   & Y^2 Y^\star -\widehat \beta Y Y^\star Y  +Y^\star X^2 -\widehat\gamma (Y Y^\star + Y^\star Y)-\widehat \rho Y^\star = \widehat\gamma^\star Y^2+\widehat\omega Y+ \widehat\eta I\,, \\
  &    (Y^\star)^2 Y -\widehat\beta Y^\star Y Y^\star  +Y (Y^\star)^2 -\widehat\gamma^\star (Y Y^\star + Y^\star Y) -\widehat \rho^\star Y= \widehat\gamma (Y^\star)^2+\widehat\omega Y^\star+ \widehat\eta^\star I\,,
\end{align}
\end{subequations}
where the parameters $\widehat\beta$, $\widehat\gamma$, $\widehat\gamma^\star$, $\widehat \rho$, $\widehat \rho^\star$, $\widehat\omega$, $\widehat\eta$ and $\widehat\eta^\star$ depend on $X$ (we recall that $X$ commutes with $Y$ and $Y^\star$). 
It would be necessary to compute the explicit dependence of the parameters in terms of $Y^\star$ (or $X$) but also to identify the relations between the four generators $X$, $Y$, $X^\star$ and  $Y^\star$. In particular, $[X^\star,Y]$ should satisfy some relations with the four previous generators. We hope to come back to this point in the future and pave the way to a proof of the above conjecture.

It would also be interesting to study the cases where the set $A_2$ is replaced by other root systems as it has been mentioned in \cite{IT}. The main other example is to replace $A_2$ by $B_2$:
\begin{align}
    B_2= \{(0,0),(1,-1), (1,0), (0,1),  (-1,1), (-1,0), (0,-1),(1,1),(-1,-1)\}\,.
\end{align}
Polynomials of ($q$-)Racah kind should appear in this case.
The higher rank factorized $A_M$-Leonard pairs introduced in Section \ref{sec:introM} could also be investigated, as well as their analogs for different root systems. In a broader viewpoint, an important question would be to understand if the approach proposed in this paper also permits to study the higher rank Leonard pairs in general and their associated multivariate polynomials. 

Another direction for generalizing the results of this paper is to consider other domains. Indeed, in this paper, we focus on the domain given by \eqref{eq:dom} but we believe that, fixing some free parameters to integer values, the domain can be restricted to some other forms. Results of this type have been obtained in \cite{IX} for Hahn multivariate polynomials.

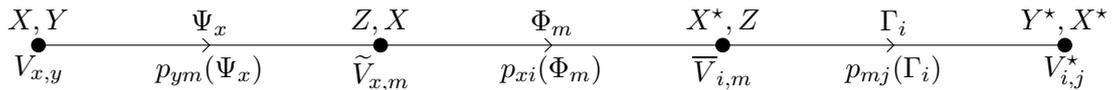
\begin{figure}[htbp]
\begin{center}
\begin{minipage}[t]{.9\linewidth}
\begin{center}
\begin{tikzpicture}[scale=0.9]
\draw[-] (-5,0)--(10,0);
\draw [fill] (-5,0) circle (0.1);
\draw (-5,0) node[above] {$X,Y$};
\draw (-5,0) node[below] {$V_{x,y}$};

\draw[-] (-2.6,0.1)--(-2.5,0)--(-2.6,-0.1);
\draw (-2.5,0) node[above] {$\Psi_{x}$};
\draw (-2.5,0) node[below] {$p_{ym}(\Psi_x)$};

\draw [fill] (0,0) circle (0.1); 
\draw (0,0) node[above] {$Z,X$};
\draw (0,0) node[below] {$\widetilde V_{x,m}$};

\draw[-] (2.4,0.1)--(2.5,0)--(2.4,-0.1);
\draw (2.5,0) node[above] {$\Phi_{m}$};
\draw (2.5,0) node[below] {$p_{xi}(\Phi_m)$};

\draw [fill] (5,0) circle (0.1);  
\draw (5,0) node[above] {$X^\star,Z$};
\draw (5,0) node[below] {$\overline{V}_{i,m}$};

\draw[-] (7.4,0.1)--(7.5,0)--(7.4,-0.1);
\draw (7.5,0) node[above] {$\Gamma_{i}$};
\draw (7.5,0) node[below] {$p_{mj}(\Gamma_{i})$};

\draw [fill] (10,0) circle (0.1);  
\draw (10,0) node[above] {$Y^\star,X^\star$};
\draw (10,0) node[below] {$V^\star_{i,j}$};

\end{tikzpicture}
\caption{Above each vertex, a couple of commuting operators is written and, below, the vectors of the basis where these operators are diagonal. At each oriented edge is associated the parameter array of the Leonard pair given by the operators above this edge.  \label{fig:trg2}}
\end{center}
\end{minipage}
\end{center}
\end{figure}

Finally, we would like to mention that the bispectrality property of the Griffiths polynomials \cite{Griff} or the Griffiths polynomials of Racah type \cite{icosi} leads to examples of rank 2 Leonard pairs which are not factorized Leonard pairs. To our knowledge, these are the only examples of polynomials which are bispectral and not of Tratnik type. To study these polynomials, it would be necessary to generalize the notion of factorized Leonard pair. In Figure \ref{fig:tr}, we provide a schematic way to understand a factorized Leonard pair. By using the same graphical representation, we believe that the structure associated to the Griffiths type polynomials is represented in Figure \ref{fig:trg2}.  This schematic interpretation has been used in \cite{icosi} as the path in the icosidodecahedron to get the Griffiths polynomials of Racah type (see the case $\partial=3$ in Figure 3 of \cite{icosi}). Let us remark that a supplementary operator, called $Z$, is necessary. We intend to start the study of this generalized factorized Leonard pair in the future.

\subsection*{Acknowledgements: }
The authors thank P.~Iliev and P.~Terwilliger for their valuable comments and LAPTh for its hospitality. They are partially supported by the international research project AAPT of the CNRS.
M.~Zaimi holds an Alexander-Graham-Bell scholarship from the Natural Sciences and Engineering Research Council of Canada (NSERC).

\appendix
\section{List of parameter arrays \label{app:pa}}

In \cite{TerClas} all the possible parameter arrays
\begin{align}
\Phi=(\theta_i,\theta^\star_i,i=0,\dots,d;\varphi_i, \phi_i,i=1,\dots,d)
\end{align} 
are listed.
We recall this list below (except for the cases Racah, $q$-Racah and Banna--Ito since we do not find $A_2$-Leonard pairs based on them) with the associated polynomials $p_{xi}(\Phi)=k_i(\Phi) u_i(\theta_x)$ ($u_i(\theta_x)$ is given by the ($q$-)hypergeometric function and $k_i(\Phi)$ is the factor in front of it) and the associated coefficients $b_i(\Phi)$ and $c_i(\Phi)$ (we recall that $a_i(\Phi)=-b_i(\Phi)-c_i(\Phi)+\theta_0$). We also provide for each parameter array  $A_2$-contiguity relations.

\subsection{List of parameter arrays of type I\label{app:pa1}}

\subsubsection{$q$-Hahn}

\paragraph{Parameter array.}

\begin{align}
&\theta_i= \theta_0+\mu(1-q^i)q^{-i}\,,
&& \theta^\star_i=\theta_0^\star+\mu^\star(1-q^i)(1-\sigma^\star q^{i+1})q^{-i} \,,\nonumber \\
&\varphi_i=\mu \mu^\star q^{1-2i}(1-q^i)(1-q^{i-d-1})(1-\rho q^i) \,,
&& \phi_i=-\mu \mu^\star q^{1-i}(1-q^i)(1-q^{i-d-1})(\rho-\sigma^\star q^{i}) \,.
\end{align}

\paragraph{Polynomial.}
\begin{align}\label{eq:polqHahn}
p_{xi}(\Phi)=  
\frac{(q\rho;q)_i (q^2 \sigma^\star ;q)_d}
{(q\rho)^i(q\sigma^\star/\rho,\sigma^\star q^{i+1};q)_i (\sigma^\star q^{2+2i};q)_{d-i}  }
\left[\begin{array}{c}d\\i\end{array}\right]_q
{}_3\phi_2 \Biggl({{q^{-i},\;\sigma^\star q^{i+1} , \; q^{-x}}\atop
{\rho q, \; q^{-d}}}\;\Bigg\vert \; q,q\Biggr)\,.
\end{align}

\paragraph{Recurrence coefficients.} 

\begin{subequations}\label{eq:recqH}
\begin{align}
   & b_i(\Phi)=\frac{\mu(1-\sigma^\star q^{i+1}) (1-q^{i-d}) (1-\rho q^{i+1})}{(1-\sigma^\star q^{2i+2})(1-\sigma^\star q^{2i+1}) }
\,,\\
   & c_i(\Phi)= \frac{ \mu q^{i-d}(1-\sigma^\star q^{d+i+1})  (1-q^{i})(\sigma^\star q^{i}-\rho) }{(1-\sigma^\star q^{2i+1})(1-\sigma^\star q^{2i})}
\,.
\end{align}
\end{subequations}

\paragraph{$A_2$-contiguity recurrence relations.}

\begin{subequations}
\begin{align}
  &\frac{q(\rho-q\sigma^\star)}{1-q^2 \sigma^\star}  p_{xi}(\Phi)   =  \frac{q^{i+1}(\rho-\sigma^\star q^{i+1})}{1-\sigma^\star q^{2i+2}}  p^{-}_{xi}(\Phi) 
  +\frac{1-\rho q^{i}}{1-\sigma^\star q^{2i}}p^{-}_{xi-1}(\Phi)  \,,\label{eq:adqh1} \\
  & (1-q^{d+1-x}) \frac{1-q\sigma^\star}{\rho-\sigma^\star} p_{xi}(\Phi)   =    \frac{q(1-q^{i+1})(1-\sigma^\star q^{i+2+d})}{1-\sigma^\star q^{2i+2}}
 p^{+}_{xi+1}(\Phi)\nonumber \\
 &\hspace{5cm} +\frac{q^{i+1}(1-q^{d+1-i})(1-\sigma^\star q^{i})}{1-\sigma^\star q^{2i}}
 p^{+}_{xi}(\Phi)\,,\label{eq:adqh2}
\end{align}
\end{subequations}
where $p^\epsilon_{xi}(\Phi)=p_{xi}(\Phi)\Big|_{ d\rightarrow d +\epsilon 1,\, \sigma^\star\rightarrow \sigma^\star/q^{\epsilon 1}}$, for $\epsilon=\pm$.

\subsubsection{Dual $q$-Hahn}

\paragraph{Parameter array.}
\begin{align}
&\theta_i= \theta_0+\mu(1-q^i)(1-\sigma q^{i+1})q^{-i}\,,
&& \theta^\star_i=\theta_0^\star+\mu^\star(1-q^i)q^{-i} \,,\\
&\varphi_i=\mu \mu^\star q^{1-2i}(1-q^i)(1-q^{i-d-1})(1-\rho q^i) \,,
&& \phi_i=\mu \mu^\star q^{d+2-2i}(1-q^i)(1-q^{i-d-1})(\sigma-\rho q^{i-d-1}) \,.\nonumber
\end{align}

\paragraph{Polynomial.}
\begin{align}\label{eq:poldqHahn}
p_{xi}(\Phi)=(-1)^i q^{\frac{i}{2}(i-3-2d)}\frac{ (q\rho;q)_i}{\sigma^i (q^{-d} \rho/\sigma;q)_i}
\left[\begin{array}{c}d\\i\end{array}\right]_q
{}_3\phi_2 \Biggl({{q^{-i}, \; q^{-x},\;\sigma q^{x+1}}\atop
{\rho q, \; q^{-d}}}\;\Bigg\vert \; q,q\Biggr)\,.
\end{align}

\paragraph{Recurrence coefficients.} 

\begin{subequations}\label{eq:recdqH}
\begin{align}
   & b_i(\Phi)=\mu(1-q^{i-d}) (1-\rho q^{i+1}) \,,\\
   & c_i(\Phi)=q\mu(1-q^i)(\sigma-\rho q^{i-d-1})
\,.
\end{align}
\end{subequations}

\paragraph{$A_2$-contiguity recurrence relations.}

\begin{subequations}
\begin{align}
  & q(\rho-\sigma q^{d})  p_{xi}(\Phi)   =  q(\rho q^{i}-\sigma q^{d})
 p^{-}_{xi}(\Phi) +(1-\rho q^{i})p^{-}_{xi-1}(\Phi) \,, \label{eq:addqh1} \\
  & \frac{(q^x-q^{d+1})(\sigma q^{d+1}-q^{-1-x})}{\rho -\sigma q^{d+1}} p_{xi}(\Phi)   =    (q^{i+1}-1) p^{+}_{xi+1}(\Phi)+(q^{d+1}-q^i)  p^{+}_{xi}(\Phi)\,,\label{eq:addqh2}
\end{align}
\end{subequations}
where $p^\epsilon_{xi}(\Phi)=p_{xi}(\Phi)\Big|_{ d\rightarrow d +\epsilon 1}$, for $\epsilon=\pm$.

\subsubsection{Dual quantum $q$-Krawtchouk}

\paragraph{Parameter array. }

\begin{align}
&\theta_i= \theta_0+\mu(q^i-1)\,,
&& \theta^\star_i=\theta_0^\star+\mu^\star(q^{-i}-1) \,,\nonumber\\
&\varphi_i= \mu \mu^\star (1-q^i)(1-q^{d-i+1})/\rho\,,
&& \phi_i=\mu \mu^\star (1-q^{-i})(1-q^{d+1-i})(1-q^{i}/\rho) \,.\label{eq:padqq}
\end{align}
We slightly modify the parameters in comparison to \cite{TerClas}: $s\to \mu/q, r \to \mu q^{d}/\rho, h^\star\to \mu^\star$. We also  change the name in comparison to \cite{TerClas} from 
 quantum $q$-Krawtchouk to dual quantum $q$-Krawtchouk
 in order to match with the usual name of the polynomial used in \cite{Koek} for example.

\paragraph{Polynomial.}
\begin{align}
p_{xi}(\Phi)=(-1)^i q^{\frac{i(i-1)}{2}} \frac{1}{(\rho q^{-i};q)_i } \left[\begin{array}{c}d\\i\end{array}\right]_q
{}_2\phi_1 \Biggl({{q^{-i}, \; q^{-x}}\atop
{\;q^{-d}}}\;\Bigg\vert \; q,\rho q^{x-d}\Biggr)\,.
\end{align}

\paragraph{Recurrence coefficients.} 

\begin{subequations}
\begin{align}
   & b_i(\Phi)=\mu(1-q^{d-i})q^{2i+1}/\rho\,,\\
   & c_i(\Phi)=\mu(1-q^i)(1- q^{i}/\rho)\,.
\end{align}
\end{subequations}

\paragraph{$A_2$-contiguity recurrence relations.}

\begin{subequations}
\begin{align}
  &  (\rho-q)p_{xi}(\Phi)   =  (\rho-q^{i+1}) p^{-}_{xi}(\Phi) + q^{i} p^{-}_{xi-1}(\Phi)  \,,\label{eq:addqqk1} \\
  & \frac{\rho q^{d+1}}{\rho-1}(1-q^{x-d-1}) p_{xi}(\Phi)   =    
  (q^{i+1}-1) p^{+}_{xi+1}(\Phi)
 +(q^{d+1}-q^{i}) p^{+}_{xi}(\Phi)\,,\label{eq:addqqk2}
\end{align}
\end{subequations}
where $p^\epsilon_{xi}(\Phi)=p_{xi}(\Phi)\Big|_{d\rightarrow d +\epsilon 1,\, \rho\rightarrow \rho q^{\epsilon 1}}$, for $\epsilon=\pm$.

\subsubsection{Quantum $q$-Krawtchouk}

We also add in the list the quantum $q$-Krawtchouk, which is not in the list of \cite{TerClas} (we recall that the 
parameter array called quantum $q$-Krawtchouk in \cite{TerClas} is called here dual quantum $q$-Krawtchouk). It is not in the usual list since it can be obtained by a transformation $q\to q^{-1}$, but it is useful here to also provide the formulas for the quantum $q$-Krawtchouk.

\paragraph{Parameter array.}

\begin{align}
&\theta_i= \theta_0+\mu(q^{-i}-1)\,,
&& \theta^\star_i=\theta_0^\star+\mu^\star(q^i-1) \,,\nonumber\\
&\varphi_i= \mu \mu^\star (1-q^i)(1-q^{d+1-i})/\rho\,,
&& \phi_i=\mu \mu^\star (1-q^{i})(1-q^{i-d-1})(1-q^{d+1-i}/\rho) \,.
\end{align}
This parameter array is obtained by the transformation \ref{eq:Phistar} of the parameter array \eqref{eq:padqq}.
 
\paragraph{Polynomial.}
\begin{align}
p_{xi}(\Phi)=(-1)^i q^{\frac{i(i+1)}{2}} \frac{1}{q^{id}(\rho q^{-d};q)_i } \left[\begin{array}{c}d\\i\end{array}\right]_q
{}_2\phi_1 \Biggl({{q^{-x}, \; q^{-i}}\atop
{\;q^{-d}}}\;\Bigg\vert \; q,\rho q^{i-d}\Biggr)\,.
\end{align}

\paragraph{Recurrence coefficients.} 

\begin{subequations}
\begin{align}
   & b_i(\Phi)=-\mu(1-q^{d-i})q^{-i}/\rho\,,\\
   & c_i(\Phi)=\mu(1-q^{-i})(1- q^{d+1-i}/\rho)\,.
\end{align}
\end{subequations}

\paragraph{$A_2$-contiguity recurrence relations.}

\begin{subequations}
\begin{align}
  &  (\rho-q^d)p_{xi}(\Phi)   =  (\rho-q^{d-i}) p^{-}_{xi}(\Phi) + q^{d+1-i} p^{-}_{xi-1}(\Phi)  \,,\label{eq:adqqk1} \\
  & \frac{\rho}{q^{d+1}-\rho}(1-q^{d+1-x}) p_{xi}(\Phi)   =    
  (1-q^{-i-1})q^{d+1} p^{+}_{xi+1}(\Phi)
 +(q^{d+1-i}-1) p^{+}_{xi}(\Phi)\,,\label{eq:adqqk2}
\end{align}
\end{subequations}
where $p^\epsilon_{xi}(\Phi)=p_{xi}(\Phi)\Big|_{d\rightarrow d +\epsilon 1 }$, for $\epsilon=\pm$.

\subsubsection{$q$-Krawtchouk}

\paragraph{Parameter array.}
\begin{align}
&\theta_i= \theta_0-\mu(1-q^{-i})\,,
&& \theta^\star_i=\theta_0^\star-\mu^\star(1-q^{-i})(1-\sigma^\star q^{i+1}) \,,\nonumber\\
&\varphi_i=\mu \mu^\star q^{1-2i}(1-q^i)(1-q^{i-d-1})\,,&& \phi_i=\mu \mu^\star \sigma^\star q(1-q^i)(1-q^{i-d-1}) \,.
\end{align}

\paragraph{Polynomial.}
\begin{align}
p_{xi}(\Phi)=(-1)^i q^{-\frac{i}{2}(3+i)} \frac{( \sigma^\star q^2;q)_d }{(\sigma^\star)^i(\sigma^\star q^{1+i};q)_i (\sigma^\star q^{2+2i};q)_{d-i}} \left[\begin{array}{c}d\\i\end{array}\right]_q
{}_3\phi_2 \Biggl({{q^{-i}, \;\sigma^\star q^{i+1},\;q^{-x}}\atop
{0,\;q^{-d}}}\;\Bigg\vert \; q,q\Biggr)\,.
\end{align}

\paragraph{Recurrence coefficients.} 

\begin{subequations}
\begin{align}
   & b_i(\Phi)=\frac{\mu(1-q^{i-d})(1-\sigma^\star q^{i+1})}{(1-\sigma^\star q^{2i+1})(1-\sigma^\star q^{2i+2})} \,,\\
   & c_i(\Phi)=\frac{\mu \sigma^\star q^{2i-d}(1-q^{i})(1-\sigma^\star q^{d+i+1})}{(1-\sigma^\star q^{2i})(1-\sigma^\star q^{2i+1})} \,.
\end{align}
\end{subequations}

\paragraph{$A_2$-contiguity recurrence relations.}

\begin{subequations}
\begin{align}
  &  p_{xi}(\Phi)   =  \frac{q^{2i}(1- q^2 \sigma^\star) }{1-\sigma^\star q^{2i+2} }   p^{-}_{xi}(\Phi) -\frac{(1- q^2 \sigma^\star)}{q^2\sigma^\star(1-\sigma^\star q^{2i})}
 p^{-}_{xi-1}(\Phi) \,, \label{eq:adqk1} \\
  & (q^{d+1-x}-1) p_{xi}(\Phi)   =  
  \frac{\sigma^\star q(1-q^{i+1})(1-\sigma^\star q^{i+d+2})}{(1-\sigma^\star q^{2i+2})(1-q\sigma^\star) }
 p^{+}_{xi+1}(\Phi)
 +\frac{\sigma^\star q(q^i-q^{d+1})(1-\sigma^\star q^i)}{(1-\sigma^\star q^{2i})(1-q \sigma^\star)}
 p^{+}_{xi}(\Phi)\,,\label{eq:adqk2}
\end{align}
\end{subequations}
where $p^\epsilon_{xi}(\Phi)=p_{xi}(\Phi)\Big|_{d\rightarrow d +\epsilon 1,\,\sigma^\star \rightarrow \sigma^\star  / q^{\epsilon 1}}$, for $\epsilon=\pm$.

\subsubsection{Affine $q$-Krawtchouk}

\paragraph{Parameter array.}
\begin{align}
&\theta_i= \theta_0-\mu(1-q^{-i})\,,&& 
\theta^\star_i=\theta_0^\star-\mu^\star(1-q^{-i}) \,,\nonumber\\
&\varphi_i=\mu \mu^\star q^{1-2i}(1-q^i)(1-q^{i-d-1})(1-\rho q^i)\,,&& \phi_i=-\mu \mu^\star \rho q^{1-i}(1-q^i)(1-q^{i-d-1}) \,.
\end{align}

\paragraph{Polynomial.}
\begin{align}
p_{xi}(\Phi)= \frac{(\rho q;q)_i}{(q\rho)^i}\left[\begin{array}{c}d\\i\end{array}\right]_q
{}_3\phi_2 \Biggl({{q^{-i}, \;0,\;q^{-x}}\atop
{\rho q,\;q^{-d}}}\;\Bigg\vert \; q,q\Biggr)\,.
\end{align}

\paragraph{Recurrence coefficients.} 

\begin{subequations}
\begin{align}
   & b_i(\Phi)=\mu(1-q^{i-d})(1-\rho q^{i+1})\,,\\
   & c_i(\Phi)=-\rho \mu q^{i-d}(1-q^i) \,.
\end{align}
\end{subequations}

\paragraph{$A_2$-contiguity recurrence relations.}

\begin{subequations}
\begin{align}
  &  p_{xi}(\Phi)   =   q^i p^{-}_{xi}(\Phi) + q^{-1}(1/\rho-q^i)p^{-}_{xi-1}(\Phi) \,, \label{eq:adaqk1} \\
  & (1-q^{d+1-x})p_{xi}(\Phi)   =  q\rho  (1-q^{i+1}) p^{+}_{xi+1}(\Phi)+ q\rho(q^i-q^{d+1}) p^{+}_{xi}(\Phi)\,,\label{eq:adaqk2}
\end{align}
\end{subequations}
where $p^\epsilon_{xi}(\Phi)=p_{xi}(\Phi)\Big|_{d\rightarrow d +\epsilon 1}$, for $\epsilon=\pm$.

\subsubsection{Dual $q$-Krawtchouk}

\paragraph{Parameter array.}
\begin{align}
&\theta_i= \theta_0+\mu(1-q^i)(1-\sigma q^{i+1})q^{-i}\,,&& \theta^\star_i=\theta_0^\star+\mu^\star(1-q^i)q^{-i} \,,\nonumber\\
&\varphi_i=\mu \mu^\star q^{1-2i}(1-q^i)(1-q^{i-d-1}) \,,&& \phi_i=\mu \mu^\star \sigma q^{d+2-2i}(1-q^i)(1-q^{i-d-1}) \,.
\end{align}

\paragraph{Polynomial.}
\begin{align}
p_{xi}(\Phi)=(-1)^i \sigma^{-i} q^{\frac{i}{2}(i-2d-3)}\left[\begin{array}{c}d\\i\end{array}\right]_q{}_3\phi_2 \Biggl({{q^{-i}, \;q^{-x},\;\sigma q^{x+1}}\atop
{0,\;q^{-d}}}\;\Bigg\vert \; q,q\Biggr)\,.
\end{align}

\paragraph{Recurrence coefficients.} 

\begin{subequations}
\begin{align}
   & b_i(\Phi)=\mu(1-q^{i-d})\,,\\
   & c_i(\Phi)=\mu \sigma q (1-q^i) \,.
\end{align}
\end{subequations}

\paragraph{$A_2$-contiguity recurrence relations.}

\begin{subequations}
\begin{align}
  &  \sigma p_{xi}(\Phi)   =   \sigma p^{-}_{xi}(\Phi) - q^{-d-1}p^{-}_{xi-1}(\Phi)  \,,\label{eq:addqk1} \\
  & (q^{-x}-q^{-d-1})(1-\sigma q^{d+x+2}) p_{xi}(\Phi)   =    \sigma q (1-q^{i+1})  p^{+}_{xi+1}(\Phi)+\sigma q^{i+1}(1-q^{d-i+1}) p^{+}_{xi}(\Phi)\,,\label{eq:addqk2}
\end{align}
\end{subequations}
where $p^\epsilon_{xi}(\Phi)=p_{xi}(\Phi)\Big|_{d\rightarrow d +\epsilon 1}$, for $\epsilon=\pm$.

\subsection{List of parameter arrays of type II\label{app:pa2}}

\subsubsection{Hahn}

\paragraph{Parameter array.}
\begin{align}
&\theta_i=\theta_0+\sigma i \,,&& \theta^\star_i=\theta^\star_0+\mu^\star i(i+1+\sigma^\star) \,,\nonumber\\
&\varphi_i=\mu^\star \sigma i(i-d-1)(i+\rho) \,,&& \phi_i=-\mu^\star \sigma i (i-d-1)(i+\sigma^\star-\rho) \,. \label{eq:parhahn}
\end{align}

\paragraph{Polynomial.}
\begin{eqnarray}\label{eq:polHahn}
p_{xi}(\Phi)=\binom{d}{i}
\frac{(\rho+1)_i(2+\sigma^\star)_d}{(\sigma^\star-\rho+1,\sigma^\star+1+i)_i(\sigma^\star+2i+2)_{d-i}} \ 
{}_3F_2 \Biggl({{-i,\;i+1+\sigma^\star, \;-x}\atop
{\rho+1,\; -d}}\;\Bigg\vert \; 1\Biggr)\,.
\end{eqnarray}

\paragraph{Recurrence coefficients.}

\begin{subequations}\label{eq:rechahn}
\begin{align}
   & b_i(\Phi)=\frac{\sigma(i-d)(i+1+\rho)(i+\sigma^\star+1)}{(2i+\sigma^\star+1)(2i+\sigma^\star+2)}\,,\\
&c_i(\Phi)=-\frac{\sigma i(i+d+\sigma^\star+1)(i+\sigma^\star-\rho)}{(2i+\sigma^\star)(2i+\sigma^\star+1)}\,.
\end{align}
\end{subequations}

\paragraph{$A_2$-contiguity recurrence relations.}

\begin{subequations}
\begin{align}
  &\frac{\sigma^\star+1-\rho}{2+\sigma^\star}  p_{xi}(\Phi)   =    \frac{i+\rho}{2i+\sigma^\star}p^{-}_{x,i-1}(\Phi)+\frac{i+1+\sigma^\star-\rho}{2i+\sigma^\star+2}p^{-}_{x,i}(\Phi)\,,\label{eq:adhahn5}\\
  &\frac{\sigma^\star+1}{\rho-\sigma^\star}(x-d-1)  p_{xi}(\Phi)   =    \frac{(d+1-i)(i+\sigma^\star)}{2i+\sigma^\star}p^{+}_{x,i}(\Phi)+\frac{(i+1)(i+\sigma^\star+2+d)}{2i+\sigma^\star+2}p^{+}_{x,i+1}(\Phi)\,,\label{eq:adhahn6}
\end{align}
\end{subequations}
where $p^\epsilon_{xi}(\Phi)=p_{xi}(\Phi)\Big|_{d\rightarrow d +\epsilon 1,\, \sigma^\star\rightarrow \sigma^\star -\epsilon 1}$, for $\epsilon=\pm$.

\subsubsection{Dual Hahn}

\paragraph{Parameter array.}
\begin{align}
&\theta_i=\theta_0+\mu i(i+1+\sigma) \,,&& \theta^\star_i=\theta_0^\star +\sigma^\star i \,,\nonumber\\
&\varphi_i=\mu \sigma^\star i(i-d-1)(i+\rho) \,,&& \phi_i=\mu \sigma^\star i(i-d-1)(i+\rho-\sigma-d-1) \,. \label{eq:pardhahn}
\end{align}

\paragraph{Polynomial.}
\begin{eqnarray}\label{eq:poldHahn}
p_{xi}(\Phi)=(-1)^i \binom{d}{i}
\frac{(\rho+1)_i}{(\rho-\sigma-d)_i}\ 
{}_3F_2 \Biggl({{-i, \;-x,\; x+1+\sigma}\atop
{\rho+1,\; -d}}\;\Bigg\vert \; 1\Biggr)\,.
\end{eqnarray}

\paragraph{Recurrence coefficients.} 

\begin{subequations} \label{eq:recdhahn}
\begin{align}
   & b_i(\Phi)=\mu(i-d)(i+1+\rho)\,,\\
&c_i(\Phi)=\mu i(i+\rho-\sigma-d-1)\,.
\end{align}
\end{subequations}

\paragraph{$A_2$-contiguity recurrence relations.}

\begin{subequations}
\begin{align}
  & (d+\sigma-\rho) p_{xi}(\Phi)   = -(i+\rho-\sigma-d) p^{-}_{xi}(\Phi)+(i+\rho)p^{-}_{xi-1}(\Phi)\,, \label{eq:addhahn3} \\
  &\frac{(x-d-1)(x+d+\sigma+2)}{d-\rho+\sigma+1}  p_{xi}(\Phi)   = (i-d-1) p^{+}_{xi}(\Phi)-(i+1)p^{+}_{xi+1}(\Phi)\,,\label{eq:addhahn4}
\end{align}
\end{subequations}
where $p^\epsilon_{xi}(\Phi)=p_{xi}(\Phi)\Big|_{d\rightarrow d +\epsilon 1}$, for $\epsilon=\pm$.

\subsubsection{Krawtchouk}

\paragraph{Parameter array.}
\begin{align}
&\theta_i=\theta_0+\sigma i \,,&& \theta^\star_i=\theta^\star_0+\sigma^\star i \,,\nonumber\\
&\varphi_i= \rho i(i-d-1)\,,&& \phi_i= (\rho-\sigma\sigma^\star)i(i-d-1)\,.
\end{align}

\paragraph{Polynomial.}
\begin{eqnarray}
p_{xi}(\Phi)=\left(\frac{\rho}{\sigma\sigma^\star-\rho}\right)^i \binom{d}{i}\
{}_2F_1 \Biggl({{-i, \;-x}\atop
{-d}}\;\Bigg\vert \; \rho^{-1}\sigma\sigma^\star\Biggr)
\,.
\end{eqnarray}

\paragraph{Recurrence coefficients.} 

\begin{subequations}\label{eq:abqkraw}
    \begin{align}
    &b_i(\Phi)=\frac{\rho}{\sigma^\star}(i-d)\,,\\
    &c_i(\Phi)=\frac{\rho-\sigma \sigma^\star}{\sigma^\star}i\,.
\end{align}
\end{subequations}

\paragraph{$A_2$-contiguity recurrence relations.}
\begin{subequations}
\begin{align}
   &p_{xi}(\Phi)=p^-_{xi}(\Phi)+\frac{\rho}{\sigma\sigma^\star-\rho}p^-_{x,i-1}(\Phi)\,, \label{eq:adkraw1}\\
   &\frac{\sigma\sigma^\star}{\rho-\sigma\sigma^\star}(x-d-1)p_{xi}(\Phi)=(d+1-i)p^+_{xi}(\Phi)+(i+1)p^+_{x,i+1}(\Phi)\label{eq:adkraw2}\,,
\end{align}
\end{subequations}
where $p^\epsilon_{xi}(\Phi)=p_{xi}(\Phi)\Big|_{d\rightarrow d +\epsilon 1}$, for $\epsilon=\pm$.

\end{document}